\input amstex
\input amsppt.sty
\magnification\magstep1

\def\ni\noindent
\def\sbs{\subset}

\def\R{\text{\bf R}}

\def\Q{\text{\bf Q}}
\def\Z{\text{\bf Z}}
\def\H{\text{\bf H}}
\def\N{\text{\bf N}}

\def\sA{\Cal A}
\def\sC{\Cal C}

\def\sP{\Cal P}

\def\sU{\Cal U}

\hoffset= 0.0in
\voffset= 0.0in
\hsize=32pc
\vsize=43pc
\baselineskip=24pt
\NoBlackBoxes
\topmatter
\author
A.N. Dranishnikov
\endauthor

\title
Asymptotic Topology
\endtitle
\abstract 
We establish some basic theorems in dimension theory and absolute extensor
theory in the coarse category of metric spaces. Some of the statements in 
this category can be translated in general topology language by applying
the Higson corona functor.
The relation of problems and results of
this `Asymptotic Topology' to Novikov and similar conjectures is discussed. 
\endabstract

\thanks The author was partially supported by NSF grant DMS-9696238.
\endthanks

\address University of Florida, Department of Mathematics, P.O.~Box~118105,
358 Little Hall, Gainesville, FL 32611-8105, USA
\endaddress


\email  dranish\@math.psu.edu
\endemail

\endtopmatter

\heading
Table of contents
\endheading

\S1 Introduction

\S2 Choice of category

\S3 Some functors in the asymptotic category

\S4 Absolute extensors

\S5 ANE, Homotopy, Dimension

\S6 Approximation by polyhedra, cohomology and dimension

\S7 Higson corona

\S8 On coarse Novikov conjectures

\S9 Principle of descent and the Higson corona

\S10 Open problems

\
\document
\head \S1 Introduction \endhead

The large scale geometry became very popular mainly thanks to Gromov.
It is a basic tool of study of discrete groups and geometric and topological
invariants of manifolds associated with the fundamental group \cite{Gr1},\cite{Ro1},
\cite{Ro2}. Objects of study in the large scale geometry are unbounded metric spaces,
such as complete open Riemannian manifolds,
their global pictures, where small (bounded) details are not taken into
 account. One investigates the properties there which are defined by taking
the limit at `infinity´.
This makes an analogy between the large scale world and the small
scale world, where one takes the limit at $0$. A significant piece
of Topology is devoted to study local properties of spaces. i.e.
the small scale world. The most developed are the theory of local connectivity,
ANR-theory, dimension theory and cohomological dimension of compact metric spaces.
The purpose of this paper is an attempt to transfer these theories to
the large scale world. Partly it was already done by other authors.
Thus the basic ideas of large scale dimension theory were introduced by
Gromov. The coarse (large scale) cohomology groups were defined by
Roe. An importance of large scale versions of concepts of local
topology was known for years. The most striking result is due to
G.Yu \cite{Y1} and it says that the Novikov higher signature conjecture holds
for geometrically finite groups with a finite large scale dimension.
We note that a discrete group $\Gamma$ has a natural metric whenever one fixes a set of
generators on it. For finitely generated groups any two such metrics
defined by means of finite sets of generators are equivalent from
the large scale point of view.

Many properties of discrete groups or universal covers of their classifying spaces
can be detected from the boundaries of the groups. So far there is no good
construction of a boundary of arbitrary group. We refer to \cite{Be} for the
axioms and some constructions.
The visual sphere at infinity is
good only in hyperbolic case. In a semi-hyperbolic case it is not a coarse 
invariant [B-R],[C-K]. Nevertheless there is the boundary of a group
 which is  coarse invariant, namely,
the Higson corona. It can be defined for any proper metric space $X$, i.e. a
metric space where the distance to a point is a proper function, and it is a 
covariant functor $\nu$ to the category of compact spaces in corresponding setting.
The problem with the Higson corona is that it is never metrizable. Thus in order
to study coarse (asymptotic) topology of nice spaces one forced to deal with
invariants of the topology of nonmetrizable compact Hausdorf spaces. Many of the technical
tools of this branch of general topology such as Marde\v si\v c factorization theorem or Shchepin's spectral theorem
are important here. The Higson corona makes a functorial connection between the large 
scale world and
an exotic small scale world which is not completely satisfactory. One of the problems is that the cohomology 
groups of Higson corona is not a homotopy invariant of coarse topology.
In the coarse topology the spaces
$\R^n$ and the $n$-dimensional hyperbolic space are homotopy equivalent whereas their Higson
coronas have different $n$-cohomologies \cite{Dr-F}.

The main motivation for studying the macroscopic topology is the Novikov higher signature conjecture.
The conjecture claims that the higher signatures of a manifold are homotopy invariant.
For aspherical manifolds it can be rephrased that the rational Pontryagin classes are homotopy invariant.
Novikov's theorem says that the rational Pontryagin classes are topological invariant [N].
Generally rational Pontryagin classes are not homotopy invariant.
The term higher signature makes its origin from the Hirzebruch signature formula:
$\langle L_k(p_1,\dots,p_k),[M^{4k}]\rangle = \sigma(M^{4k})$ where $L_k(p_1,\dots,p_k)$ is the
$k$-th Hirzebruch polynomial in Pontryagin classes and $[M^{4k}]$ is the fundamental class.
For a simply connected manifold $M$ the signature $\sigma(M)$ is the only homotopy invariant
which can be written in terms of Pontryagin classes. In the presence of the fundamental
group $\Gamma=\pi_1(M)$ there are additional possibilities. In that case besides the
fundamental class $[M]$ there are other 'special' homology classes which comes from the
homology of the group $\Gamma$. Integration of the Hirzebruch polynomials on these classes
defines the {\it higher signatures} of a manifold.

The restriction to the case of aspherical manifolds $M$ is not that strong as it seems.
It covers the Novikov conjecture for manifolds with all geometrically finite fundamental groups.
We recall that a group $\Gamma$ is {\it geometrically finite} if $B\Gamma=K(\Gamma,1)$ is
a finite complex. In that case the universal cover $E\Gamma$ is contractible and coarsely equivalent to
$\Gamma$.

Another serious conjecture about aspherical manifolds is the Gromov-Lawson conjecture [G-L]:
{\it An aspherical manifold cannot carry a metric of a positive scalar curvature}.
The connection between this and the Novikov conjecture is discussed in [Ros], [F-R-R].

Let $X$ be the universal cover of an aspherical manifold $M^n$ supplied with a metric lifted from $M^n$.
Then the above venerable conjectures can be reduced to some large scale problems about
the metric space $X$. As a topological space $X$ is a contractible $n$-manifold. Hence
$X\times\R$ is homeomorphic to $\R^{n+1}$. Without loss of generality we may assume in the above
problems that $X=\R^n$. All problems are concentrated in the metric on $X$.
In \S8 we consider several conjectures about $X$, the validity of which would imply the
Novikov and the Gromov-Lawson conjectures.

It turns out that several statements of general topological nature formulated in
macroscopic world imply the Novikov (Gromov-Lawson) conjecture. We discuss here mainly dimensional
theoretic statements. As it was mentioned, Yu proved [Y1] that if the macroscopic dimension of
a geometrically finite group with a word metric on it is finite then the Novikov and the
Gromov-Lawson conjecture holds for manifolds with the fundamental group $\Gamma$.
In this paper we define the notion of macroscopic cohomological dimension $as\dim_{\Z}$
and show that $as\dim_{\Z}\Gamma<\infty$ for all geometrically finite groups.
This reduces the Novikov conjecture to the old Alexandroff problem about an equivalence of
two concepts of dimension: geometrical (Lebesgue's dimension) and algebraic (cohomological dimension)
considered in the macroscopic world. The Alexandroff Problem in classical
dimension theory was resolved by a counterexample [Dr]. This counterexample was also 
transfered to the macroscopic world in [D-F-W]
where we defined a uniformly contractible Riemannian metric on $\R^8$
in such a way that the asymptotic dimension is infinite and the macroscopic cohomological
dimension is finite. Nevertheless in the large scale topology
the problem is still open in the case when a metric space is a group with
word metric on it. With respect to that it would be useful to investigate the situation
in classical case. For compact topological groups the equivalence of these two dimensions
easy follows from their presentation as inverse limits of Lie groups.
It is unclear whether compact groups in the microscopic world are proper analogs of
discrete groups in the macroscopic world. It seems that the Alexandroff Problem for compact
H-spaces is more relevant to the large scale world.

The body of the paper consists of \S 2-9. In \S 10 we present a list of open problems.

In \S 2 we define a category suitable for the macroscopic topology.
We called it the Asymptotic category.
Note that our category differs slightly from Roe's Coarse category.
This difference is stipulated by our desire to build an analog of the ANR theory in
the macroscopic world.

In \S 3 we define some basic constructions in the asymptotic category.

In \S 4 the notion of AE (absolute extensor) is defined and simple examples of AE objects are
presented.

In \S 5 we define ANE spaces and  a notion of homotopy between two morphisms.
Also we discuss different approaches to definition of dimension in the
asymptotic category. In particular we define a macroscopic dimension in terms of
extension of maps to objects analogous to spheres.

In \S 6 we define the coarse cohomology by means of anti-\v Cech approximations by
polyhedra and then we define an asymptotic cohomological dimension.

In \S 7 we consider the Higson corona $\nu X$ and compare various forms of the 
asymptotic dimension of a metric space $X$ with the covering dimension of 
the Higson corona $\nu X$. The main result here is that the categorical definition
of asymptotic dimension of $X$ agrees with the covering dimension of 
$\nu X$.

In \S 8  we prove some versions of coarse Novikov conjectures under some dimensional
conditions on groups. Our results here are extensions of Yu's [Y1],[Y2].
In particular we prove  the Coarse Baum-Connes conjecture for groups with the slow 
dimension growth.

In \S 9 we compare different coarse reduction of the rational Novikov
conjecture in the language of the Higson corona.

\head \S2 Choice of category \endhead

A metric space $(X,d)$ is called {\it proper} if every closed ball
$B_r(x_0)=\{x\in X\mid d(x,x_0)\le r\}$ is compact.
A map $f:X\to Y$ is called {\it proper} if $f^{-1}(C)$ is compact
for every compact set $C\subset Y$. Note that a proper map between
proper metric spaces is always continuous.

A map $f:(X,d_X)\to (Y,d_Y)$ is called asymptotically Lipschitz if there are
numbers $\lambda$ and $s$ such that $d_Y(f(x),f(x'))\le \lambda d_X(x,x')+s$.

An {\it asymptotic category} $\sA$ consists of proper metric spaces and
proper asymptotically Lipschitz maps.

We consider the following refinement $\tilde\sA$ of this category:
Let $x_0\in X$ be a base point, then the norm $\|x\|$ of $x\in X$ is the
distance $d_X(x,x_0)$. Then the norm $\|f\|$ of a map $f:X\to Y$ is
$\lim_{x\to\infty}\frac{\|f(x)\|}{\|x\|}$. It is clear that the
norm $\|f\|$ does not depend on choice of base points.
A map with a nonzero norm is always proper. In the category $\tilde\sA$ we 
consider morphisms of $\sA$ with nonzero norms.
For any such map $f:X\to Y$ there are constants $c$ and $b$ such that
$\|f(x)\|\ge c\|x\|-b$.

EXAMPLE 1. The Euclidean $n$-space $\R^n$, the half $n$-space $\R^n_+$ with
the standard metric are examples of objects in $\sA$. 

Following Gromov we are going to call Lipschitz maps with the Lipschitz 
constant one as {\it short maps}.
A distance to a given
point is a typical morphism in $\tilde\sA$, $d(-,x_0):X\to\R_+$ on a proper metric space $(X,d)$
to $\R_+$. It is a short map. 

\proclaim{Proposition 2.1} If $f:X\to Y$ is asymptotically Lipschitz with the
Lipschitz constant $\lambda$, then for any $\alpha>1$ there is $R>0$ such that
$d_Y(f(x),f(x'))<\alpha\lambda d_X(x,x')$ provided $d(x,x')\ge R$.
\endproclaim

Isomorphisms $f:X\to Y$ in the category $\sA$ are homeomorphisms with $f$ and
$f^{-1}$ asymptotically Lipschitz. Isomorphisms in $\tilde\sA$ are isomorphisms
in $\sA$ with nonzero norm.

EXAMPLE 2. All Banach $n$-dimensional spaces are isomorphic to $\R^n$ in $\tilde\sA$
(and hence in $\sA$).

There is a bigger category $\bar\sA$ which is also of great importance. The
objects of category $\bar\sA$ are the same as in $\sA$ (actually, one can take 
all metric spaces) and morphisms are asymptotically 
Lipschitz and {\it coarsely proper} maps. The last means that $f^{-1}(C)$ is
bounded for every bounded set $C\subset Y$. In this case $f$ is not necessarily
continuous. One can define a ~ refinement of $\bar\sA$.
Two morphisms in $\bar\sA$ are called {\it coarse equivalent}
if they are within a finite distance i.e. there is $c>0$ such that
$d_Y(f(x),g(x))<c$ for all $x\in X$. A morphism $f:X\to Y$ is called
a {\it coarse isomorphism} if there is a morphism $g:Y\to X$ s.t.
$f\circ g$ and $g\circ f$ are equivalent to $1_X$ and $1_Y$ respectively.
Metric spaces $X$ and $Y$ are coarse isomorphic (quasi-isometric) if there is
a coarse isomorphism $f:X\to Y$.
The quotient category $\sC=\bar\sA/\sim$ is called the {\it coarse category}.

We recall that the Gromov-Hausdorff distance between metric spaces $d_{GH}(X,Y)$
is the infimum of distances between subsets $i(X),j(Y)\subset Z$ for all
possible isometric imbeddings $i,j$ for all possible metric spaces $Z$.
Very often $d_{GH}(X,Y)=\infty$.
\proclaim{Proposition 2.2}
If $d_{GH}(X,Y)<\infty$ then $X$ and $Y$ are coarsely isomorphic.
\endproclaim
Note that a map which is only a finite distance apart from asymptotically 
Lipschitz map is asymptotically Lipschitz itself.

DEFINITION. A metric space $X$ is called {\it uniformly contractible},
we denote $X\in UC$, if there is a function $S:[0,\infty)\to\R$ such that
every ball $B_r(x)$ of radius $r$ centered at $x$ can be contracted 
to a point in the ball $B_{S(r)}(x)$.
\proclaim{Proposition 2.3}
Let $Z,X$ be objects of $\sA$ and let $\dim Z<\infty$ and $X\in UC$. 
Then every morphism $f:Z\to X$ in $\bar\sA$ is coarsely isomorphic to
a morphism in $\sA$.
\endproclaim 
\demo{Proof}
Let $\dim Z=n$, let $\lambda$ and $s$ be the constants of $f$ from the
definition of asymptotically Lipschitz condition and let $S$ be a function from
the definition of $UC$ property of $X$. 
 Let $\sU$ be an open covering of $Z$ of order $\le n+1$ by sets
of diameter $\le 1$. Let $\phi:Z\to N$ be a projection to the nerve of
the covering $\sU$. We define a map $g_0:N^{(0)}\to X$ on $0$-dimensional skeleton
of $N$ by the rule: $g_0(u)=f(x_u)$ for some $x_u\in\phi^{-1}(u)\subset U$.
The $UC$ property allows to extend $g_0$ to a map $g:N\to X$ in such way that
the image of a $k$-dimensional simplex $g(\sigma^k)$ is contained in the
$S_k$-neighborhood of $g_0((\sigma^k)^{(0)})$, where $S_k=S(2S(\ldots S(\lambda+s)$
is a number obtained after $k$ iterations. Then for given $z\in Z$ we have
$d(f(z),g\phi(z))\le d(f(z),g_0(u))+d(g_0(u),g\phi(z))$. We take $u$ such
that $\phi(z)$ and $u$ are in the same simplex. Then
$d(f(z),g\phi(z))\le d(f(z),f(x_u))+S_n\le 1+S_n=c$.\qed
\enddemo
REMARK. In the above proof instead of $\dim Z<\infty$ we can use a weaker 
condition: $Z$ admits a cover of finite order by uniformly bounded sets.

In [Roe2] the coarse category is defined with slightly different morphisms.
A map $f:X\to Y$ is a coarse map if it is coarse proper and coarsely uniform i.e. there is a
positive function $S:[0,\infty)\to\R$ such that

$d_X(x,x')\le r\Rightarrow d_Y(f(x),f(x'))\le S(r)d_X(x,x')$.

It is clear that this property is weaker than asymptotically Lipschitz condition.

We recall that a metric space $(X,d)$ is called a {\it geodesic metric space}
if for every two points $x,y\in X$ there is an isometric imbedding
$j:[0,d(x,y)]\to X$ with $j(0)=x$ and $j(d(x,y))=y$.
\proclaim{Proposition 2.4}
Let $X$ be a geodesic metric space, then every coarse map $f:X\to Y$ is
asymptotically Lipschitz.
\endproclaim
\demo{Proof}
We take $s=S(1)$ and $\lambda=S(1)$. If $d(x,x')=n+\alpha$, $0\le\alpha<1$,
then $d_Y(f(x),f(x')\le d_Y(f(x),f(j(1)))+\dots+d_Y(f(j(n-1)),f(j(n)))+d_Y(f(j(n)),f(x'))
\le S(1)n+S(1)\le \lambda n+s$. Here $j:[0,d(x,x')]\to X$ is a geodesic segment
joining $x$ and $x'$.

\enddemo
Nevertheless there are situations when Roe's morphisms are more appropriate.
For example, they give a richer embedding theory.
A morphism $j:Y\to X$ is called an {\it imbedding} in $\sA$ if $j$ is an injection,
$j(Y)$ is a closed subset and $j^{-1}:(j(Y),d_X\mid_{j(Y)})\to Y$ is a morphism 
in $\sA$, i.e. asymptotically Lipschitz. Thus, a bending  a line into parabola with
the induced metric from the plane is not an embedding of a line into $\R^2$ in $\sA$ 
and it is an embedding in Roe' sense.

The main flaw of coarsely uniform maps is that there is no good extension theory
for them (see Remark 2 in \S4). We prefer asymptotically Lipschitz maps in the asymptotic category
because coarsely uniform maps makes weaker the analogy between proper metric spaces and
compacta. Anyway, Proposition 2.4 shows that our category and Roe's are 
very close.

\head \S3 Some functors in the asymptotic category \endhead

{\bf 1. Product.} The {\it Cartesian product} $X\times Y$ of two spaces in $\sA$ is well defined.
The metric on the product can be taken as follows $d((x_1,y_1),(x_2,y_2))=
d_X(x_1,x_2)+d_Y(y_1,y_2)$ which is equivalent to the `euclidean' metric
$\sqrt{d_X(x_1,x_2)^2+d_Y(y_1,y_2)^2}$. The problem with this product is that
it is not categorical. The projections onto the factors are not morphisms in $\sA$.
Here we define an {\it asymptotic product} $X\tilde\times Y$ as pull-back in
the topological category in the following diagram:
$$
\CD
X\tilde\times Y @>f>> Y\\
@VgVV @Vd_YVV\\
X @>d_X>> \R_+\\
\endCD
$$
Here $d_X(x)=d_X(x,x_0)$ and $x_0\in X$ is a base point. The metric on
$X\tilde\times Y$ is taken from the product $X\tilde\times Y\subset X\times Y$.
This definition is good for connected spaces. Generally this operation could give
 an empty space. 
\proclaim{Proposition 3.1} Let $X$ and $Y$ be geodesic metric spaces and
let $Z(x_0,y_0)$ denote the asymptotical product $X\tilde\times Y$ defined
by means of the base points $x_0\in X$ and $y_0\in Y$. Then for any points
$x'_0\in X$ and $y_0'\in Y$ the spaces $Z(x_0,y_0)$ and $Z(x_0',y_0')$ are in
a finite Gromov-Hausdorf distance.
\endproclaim
\demo{Proof}
Let $(x,y)\in Z(x_0,y_0)$. Then by the definition $d_X(x,x_0)=d_Y(y,y_0)$.
Let $d=\min\{d_X(x,x_0), d_X(x,x_0'),d_Y(y,y_0')\}$.
Consider geodesics $[x'_0,x]$ in $X$ and $[y_0',y]$ in $Y$ and take points
$x'\in[x_0',x]$ and $y'\in[y_0',y]$ with $d_X(x'_0,x')=d=d_Y(y'_0,y')$.
Then $(x',y')\in Z(x_0',y_0')$. Note that $d_{X\times Y}((x,y),(x',y'))=
d_X(x,x')+d_Y(y,y')=d_X(x,x'_0)-d+d_Y(y,y'_0)-d\le
|d_X(x_0',x)-d_X(x,x_0)|+|d_Y(y_0',y)-d_Y(y,y_0)|\le
d_X(x_0,x_0')+d_Y(y_0,y_0')=R$. Therefore $(x,y)$ lies in the $R$-neigborhood of
$Z(x_0',y_0')$. Hence $Z(x_0,y_0)$ lies in the $R$-neighborhood of
$Z(x_0',y_0')$. Similarly, $Z(x_0',y_0')$ lies in the $R$-neighborhood of
$Z(x_0,y_0)$.\qed
\enddemo
 Thus,by Proposition 2.2 the class of
 $X\tilde\times Y$ in the coarse category does not depend on choice 
 of base point.
 
As it will be seen latter, the half line $\R_+$ plays the role of point in
the category $\sA$. Hence the above definition is quite logical.
Note that $X\tilde\times\R_+=X$.

Since the product with unit interval turns into the identity functor after
coarsening, we need a modified notion of the product with interval. 
Again as it will be shown latter, the half plane $\R^2_+$ is an analog
of the unit interval in $\sA$. Then $X\tilde\times\R^2_+$ substitutes for
the product with unit interval.
Let $x_0\in X$ be a base point.
We define a map $J:X\tilde\times\R^2_+\to X\times\R$ by the formula:
$J(x,(s,t))=(x,s)$ where $t\in\R_+$ and $s\in\R$.
\proclaim{Proposition 3.2}
$J$ is an imbedding in the category $\sA$.
\endproclaim
\demo{Proof}
First, the image $im J$ consists of the region in $X\times\R$ between
the graphs of the functions $d_X$ and $-d_X$. Show that 
$J:X\tilde\times\R^2_+\to im J$ is
bi-Lipschitz. Note that $d_{X\times\R}(J(x,(s,t)),J(x',(s',t')))=
d_X(x,x')+|s-s'|\le 
d_X(x,x')+|s-s'|+|t-t'|=d_{X\times\R^2_+}((x,(s,t)),(x',(s',t')))$.
On the other hand, since $\|x\|=t+|s|$ and $\|x'\|=t'+|s'|$, we have
 
$2d_{X\times\R_+}(J(x,(s,t)),J(x',(s',t')))=2d_X(x,x')+2|s-s'|\ge
d_X(x,x')+|\|x\|-\|x'\| |+||s|-|s'|| +|s-s'|\ge d_X(x,x')+|t-t'|+|s-s'|=
d_{X\times\R^2_+}((x,(s,t)),(x',(s',t')))$.\qed
\enddemo
We define two maps $i_{\pm}:X\to X\tilde\times\R^2_+$ by the formula
$i_{\pm}(x)=J^{-1}(x,\pm\|x\|)$.

\proclaim{Proposition 3.3}
The maps $i_{\pm}$ are imbeddings in the category $\tilde\sA$
(and hence in $\sA$).
\endproclaim
If spaces $X$ and $Y$ have base points, one can define the {\it wedge} $X\vee Y$ by setting
$d(x,y)=d_X(x,x_0)+d_Y(y_0,y)$ for $x\in X$ and $y\in Y$. 

{\bf 2. Quotient space.}
Let $f:X\to Y$ be a proper map of a metric space. We define a metric $d_Y$
({\it quotient metric}) on $Y$ such that $f$ will be a morphism.
First we define a function $s_Y:Y\times Y\to\R_+$ by the
formula 

$s_Y(y,y')=d_X(f^{-1}(y),f^{-1}(y'))=\inf\{d_X(x,x')\mid x\in f^{-1}(y),
x'\in f^{-1}(y')\}$.

Since a proper map is always surjective, the function $s_Y$ is well defined.
Now we define a metric $d_Y$ as the intrinsic metric generated by $s_Y$:

$d_Y(y,y')=\inf\{\Sigma_{i=1}^ns_Y(y_i,y_{i+1})\mid n\in\Bbb N; y_1,\dots y_n\in Y,
y_1=y, y_n=y'\}$.

Then $d_Y(f(x),f(x'))\le s_Y(f(x),f(x'))\le d_X(x,x')$ and hence $f$ is
a short map.

{\bf 3. Cone.} In this section we define the {\it cone} over a proper metric 
space X. We identify $X\tilde\times\R^2_+$ with its image under the imbedding $J$.
Then on $i_+(X)=\{(x,\|x\|)\in X\times\R_+\}$ we consider a function $f'$ defined
by the formula $f'((x,\|x\|))=\|x\|$. We extend this map to a proper map
$f:X\tilde\times\R^2_+\to Y$ by adding singletons to the decomposition
generated by $(f')^{-1}$. Then by the definition
the quotient space $Y$ with the quotient metric $d_Y$ is the cone $CX$ over $X$.

\proclaim{Lemma 3.4}
The cone $C(\R^n)$ is isomorphic to $\R_+^{n+1}$ in $\tilde\sA$
 (and hence in $\sA$).
\endproclaim
\demo{Proof}
Let $q:J(\R^n\tilde\times\R^2_+)\to C\R^n$ be the quotient map.
We define a map $T:C(\R^n)\to\R^n\times\R_+$ by the formula

$T(q(x,t))=(\frac{\|x\|-t}{\|x\|}x,t)$.

Then the inverse transform is given by the formula:

$T^{-1}(y,t)=q((\|y\|+t)\frac{y}{\|y\|},t)$.

First we show that the inverse map is Lipschitz. Denote by $d$ the metric
$d_{\R^n\times\R_+}$. Then

\

$d_{C(\R^n)}(T^{-1}(y,t),T^{-1}(y',t'))=d_{C(\R^n)}(q((\|y\|+t)\frac{y}{\|y\|},t),
q((\|y'\|+t')\frac{y'}{\|y'\|},t'))\le $

\

$\le d(((\|y\|+t)\frac{y}{\|y\|},t),
((\|y'\|+t')\frac{y'}{\|y'\|},t'))=\|\|y\|+t)\frac{y}{\|y\|}-
\|y'\|+t')\frac{y'}{\|y'\|}\|+|t-t'|$

we may assume that $\|y'\|\le\|y\|$, then we continue:

$\le (1+\frac{t}{\|y\|})\|(y-y')\|+\|y'\||\frac{t}{\|y\|}-\frac{t'}{\|y'\|}|+|t-t'|
\le $

\ 
$2\|y-y'\|+|t-t'|+\frac{\|y'\|}{\|y\|}|t-t'|+\frac{t'}{\|y\|}\|y-y'\|\le
3\|y-y'\|+3|t-t'|=3d((y,t),(y',t'))$.

\

Let $(x,t)\in J(\R^n\tilde\times\R^2_+)$. 
Note that $d((x,t),i_+(\R^n))=\|x\|-t$. 
It is clear that 

\

$d_{C(\R^n)}(q(x,t),q(x',t'))\ge\min\{d((x,t),(x',t')),
d((x,t),i_+(\R^n))+d((x't'),i_+(\R^n))\}$.

\

We may assume that $\|x'\|\ge\|x\|$. If $t'\le\|x\|$, then $t'-t\le\|x\|-t$.
Since always $t-t'\le \|x'\|-t'$, we have

$2(d((x,t),i_+(\R^n))+d((x't'),i_+(\R^n)))=2(\|x\|-t+\|x'\|-t')\ge$

\

$\|x\|-t+\|x'\|-t'+|t-t'|\ge d(T(q(x,t)),T(q(x',t')))$.

\

On the other hand, the above argument shows that
$d(T(q(x,t)),T(q(x',t')))\le 3\|x-x'\|+3\|t-t'\|=3d((x,t),(x',t'))$.

\

Hence, $d(T(q(x,t)),T(q(x',t')))\le 3d_{C(\R^n)}(q(x,t),q(x',t'))$,
provided $t'\le\|x\|$.

Now we assume that $t'\ge\|x\|$. It easy to check that for any two
points $(z_1,t_1)$ and $(z_2,t_2)$ lying in $i_+(\R^n)$ the sum of
distances 

\

$d_{C(\R^n)}(q(x,t),q(z_1,t_1))+d_{C(\R^n)}(q(z_2,t_2),q(z_1,t_1))
+d_{C(\R^n)}(q(x',t'),q(z_2,t_2))$ 

\

greater or equal than $|\|x'\|-\|x\||+|t'-t|$.

\

Therefore, $d_{C(\R^n)}(q(x,t),q(x',t'))\ge |\|x'\|-\|x\||+|t'-t|=
\|x'\|-\|x\|+t'-t$. 

\

Since $2(t'-t)+\|x'\|-\|x\|\ge \|x\|-t+\|x'\|-t'+t'-t$,
it follows that 

\

$2d_{C(\R^n)}(q(x,t),q(x',t'))\ge d(Tq(x,t), Tq(x',t'))$.
\qed
\enddemo

One can similarly define the suspension $\Sigma X$.

{\bf 4. Probability measures.} First we recall that probability measures form a functor $P:\sA\to\sA$.
The set of all probability measures with compact supports $P(X)$ on a metric
space $(X,\rho)$ can be given a metric space structure by means
the Kantorovich-Rubinshtein metric $\bar\rho(\mu_1,\mu_2)=\sup\{|\int fd\mu_1-\int
fd\mu_2|\mid f\in S(X)\}$ where $S(X)$ is the set of all short real valued functions.
i.e. Lipschitz functions with Lipschitz constant one. Note that $X$ is 
isometrically imbedded in $P(X)$ by means of Dirac measures $\delta_x$.
For every $n$ there is a subfunctor $P_n:\sA\to\sA$ of probability measures
supported at most by $n$ points. 
Let $(X,d)$ be a proper metric space with a base point $x_0$. 
The half line space $\R_+$ has the natural base point $\{0\}$.  
One can define the {\it join}
product $X\ast Y$ of two spaces with base points $X$ and $Y$ as a subfunctor of
$P_2(X\vee Y)$. The natural question here whether one can interpret the cone
$CX$ as the join of $X$ and  $\R_+$.

For a proper metric space $X$ we denote by $L(X)$ the space of all 
short continuous functions on $X$ with the $\sup$ norm. Note that $X$ 
can be isometrically imbedded in $L(X)$.

\head \S4 Absolute extensors \endhead

For every category where the notion of subobject is defined, one can
consider absolute extensors. An object $Y$ in Category $\sC$ is an {\it absolute 
extensor} in $\sC$, $Y\in AE(\sC)$, if any other object $X$ and subobject 
$A\subset X$ and
a morphism $\phi:A\to Y$ there is an extension $\bar\phi:X\to Y$.

A subspace of $(X,d)$ in the category $\sA$ is a closed subset with the induced
metric $(Z,d|_Z)$, $Z\subset X$. What are $AE(\sA)$ and $AE(\tilde\sA)$? The one point space is
not an absolute extensor because unbounded spaces cannot have proper maps to
a point.

Let $\R^n_+=\{(x_1,\dots x_n)\in\R^n\mid x_n\ge 0\}$ be a half space and 
denote $\R_+=\R_+^1$.
\proclaim{Theorem 4.1}
$\R_+\in AE(\sA)$ and $\R_+\in AE(\tilde\sA)$.
\endproclaim
\demo{Proof}
First we consider the case $\tilde\sA$. Let $A\subset X$ be a closed subset
and let $\phi:A\to\R_+$ be a morphism with the coarse Lipschitz constants 
$\lambda$ and $s$. Let $\|\phi(x)\|\ge c\|x\|-b$. We take $c\le\lambda$. 
We choose 0 as the base point in $\R_+$. Using the idea of [CGM] we define
not necessarily continuous extension $\phi'$ of $\phi$ by transfinite induction.
Enumerate points of $X\setminus A$. Assume that $\phi'$ is already defined
on a set $B_{\alpha}=A\cup(\cup_{\beta<\alpha}\{x_{\beta}\})$ with the same constants
$\lambda$, $s$, $c$ and $b$. In order to maintain the asymptotical Lipschitz 
inequality we should take the point $x_{\alpha}$ to a point lying in all
intervals $I_x=[\phi'(x)-\lambda d(x,x_{\alpha})-s,\phi'(x)+\lambda d(x,x_{\alpha})+s)]$.
for $x\in B_{\alpha}$. Since for $\phi'(x')\le\phi'(x)$ the inequality
$\phi'(x)-\phi'(x')\le \lambda d(x,x')+s$ and the triangle inequality 
for $x$, $x'$ and $x_{\alpha}$ imply that
$\phi'(x)-\lambda d(x,x_{\alpha})-s\le \phi'(x')+\lambda d(x',x_{\alpha})+s$,
it follows that $I_x\cap I_{x'}\ne\emptyset$. In order to have inequality
$\|\phi'(x_{\alpha})\|=\phi'(x_{\alpha})\ge c\|x_{\alpha}\|-b$ we should take the point $x_{\alpha}$
to the interval $J_{\alpha}=[c\|x_{\alpha}\|-b,\infty)$. Note that
$\phi'(x)+\lambda d(x,x_{\alpha})+s\ge c\|x\|-b +\lambda d(x,x_{\alpha})+s\ge
c(\|x\|+d(x,x_{\alpha}))-b+s\ge c\|x_{\alpha}\|-b+s\ge c\|x_{\alpha}\|-b$.
Hence all intervals $I_x\cap J_{\alpha}$ are having pairwise nonempty
intersections. By Helly Theorem the intersection $I=(\cap I_x)\cap J_{\alpha}$
is nonempty. Define $\phi'(x_{\alpha})\in I$. 

By a relative version of 
Proposition 2.3 there is a continuous extension $\bar\phi$ which is in a finite
distance with $\phi'$.

We cannot apply the above argument in the case of $\sA$, since we do not have
a control on properness of extending map. Here we use different argument.
 
Let $A\subset X$ be a closed subset with the induced metric and let
$\phi:A\to\R_+$ be a proper asymptotically Lipschitz map with constants $\lambda$
and $s$. Let $R$ be as in Proposition 2.1 for $\alpha=2$. We set $m=\lambda R+s$.
For all $i$ we define $A_i=\phi^{-1}([0,mi])$ and $B_i=A\setminus A_i$.
For any positive function $\xi:E\to\R_+$, $E\subset X$, we denote by
$N_{\xi}(E)=\cup_{y\in E}B_{\xi(y)}(y)$, the $\xi$-neighborhood of $E$.
An open $\xi$-neighborhood of $E$ we denote as $ON_{\xi}(E)$. If $\xi$ is
a constant, say $\epsilon$, then these are ordinary $\epsilon$-neighborhoods.
We define a function $\xi_i$ by the formula $\xi_i(y)=\frac{1}{2\lambda}(\phi(y)-mi)$.
It is easy to check that $\xi_i>0$ on $B_{i+1}$. Denote $O_i=ON_{\xi_i}(B_{i+1})$.
By induction we construct a family $\{C_i\}$ of a compact subsets of $X$ with
the properties:
\roster
\item{} $C_i\cap A=A_i$ and $Int(C_i)\cap A=Int_A(A_i)$,
\item{} $N_{\frac{m}{2\lambda}}(C_i)\subset C_{i+1}$,
\item{} $C_i\cap O_i=\emptyset$,
\item{} $N_{\frac{m}{2\lambda}}(C_i)\cap B_{i+1}=\emptyset$
\endroster
Conditions (3),(4) are technical, they are needed only for a smooth work of
induction.

We start with $C_0=A_0$. Now assume that we have constructed
$C_0,\dots, C_i$ satisfying (1)-(4). Since $\phi$ is a proper map, 
$A_{i+1}$ is a compact set. There is a compact set $C'$ satisfying (1).
Note that $C=C'\cup N_{\frac{m}{2\lambda}}(C_i)$ satisfies (1) as well.
Indeed, $C\cap A=(C'\cap A)\cup(N_{\frac{m}{2\lambda}}(C_i)\cap A)=A_{i+1}\cup
(N_{\frac{m}{2\lambda}}(C_i)\cap A_{i+1})\cup(N_{\frac{m}{2\lambda}}(C_i)\cap B_{i+1})=
A_{i+1}$ by virtue of (4). We define $C_{i+1}=C\setminus O_{i+1}$.

To check (1) it suffices to show that $Cl(O_{i+1})\cap A_{i+1}=\emptyset$.
If $x\in Cl(O_{i+1})\cap A_{i+1}$, then there is $y\in Cl(B_{i+2})$ with
$d(x,y)\le\xi_{i+1}(y)$. Since $\phi(x)\le m(i+1)$ and $\phi(y)\ge m(i+2)$,
$\phi(y)-\phi(x)\ge m$. Hence $\lambda d(x,y)+s\ge m=\lambda R+s$, i.e.
$d(x,y)\ge R$. Therefore, $\phi(y)-\phi(x)< 2\lambda d(x,y)$. On the other
hand, $\phi(y)-\phi(x)\ge \phi(y)-m(i+1)=2\lambda\xi_{i+1}(y)\ge 2\lambda 
d(x,y)$. The contradiction implies that $Cl(O_{i+1}\cap A_{i+1}=\emptyset$.

To check (2), it suffices to show that $N_{\frac{m}{2\lambda}}(C_i)\cap 
O_{i+1}=\emptyset$. Assume the contrary, there is 
$x\in N_{\frac{m}{2\lambda}}(C_i)\cap O_{i+1}$. Then there are points
$z\in C_i$ and $y\in B_{i+2}$ such that $d(x,z)\le m/2\lambda$ and
$d(x,y)<\xi_{i+1}(y)$. Then by the triangle inequality
$d(z,y)<\xi_{i+1}(y)+m/2\lambda=\frac{1}{2\lambda}(\phi(y)-mi)=\xi_i(y)$. Hence,
$z\in O_i$. Thus, $z\in C_i\cap O_i$ which contradicts with (3).

The condition (3) is satisfied automatically. Check the condition (4).
Assume the contrary: there is $x\in N_{\frac{m}{2\lambda}}(C_{i+1})\cap B_{i+2}$.
Then there is $y\in C_{i+1}$ with $d(x,y)\le m/2\lambda$. On the other hand,
by the condition (3) for $C_{i+1}$ we have $d(x,y)\ge\xi_{i+1}(x)=
1/2\lambda(\phi(x)-m(i+1))$. These inequalities imply $m(i+2)\ge \phi(x)$.
So, $x$ cannot be in $B_{i+2}$. Contradiction.

We note that the condition (2) implies that $\cup C_i=X$.
We define $\phi(\partial C_i)=mi$. Because of the condition (1) this is
a continuous extension of $\phi$ over $A'=A\cup(\cup_{i=0}^{\infty}\partial C_i)$.
Consider the sets $D_i=C_i\setminus Int(C_{i-1})$. Extend the map
$\phi|_{D_i\cap A'}:D_i\cap A'\to[m(i-1),mi]$ to an arbitrary continuous map
$\bar\phi_i:D_i\to[m(i-1),mi]$. The union $\bar\phi=\cup\bar\phi_i$ is a proper
continuous map $\bar\phi:X\to\R_+$. Show that $\bar\phi$ is asymptotically
Lipschitz. Take two points $x,y\in X$. Take the minimal number $i$ such that
$x\in C_i$. Similarly let $j$ be the minimal number such that $y\in C_j$.
With out loss of generality we may assume that $j=i+k$. Then
$|\bar\phi(y)-\bar\phi(x)|\le m(k+1)$. By the condition (2) we have 
$d(x,y)\ge (k-1)m/2\lambda$. Therefore, $k+1\le (2\lambda/m)d(x,y)+2$.
Hence, $|\bar\phi(y)-\bar\phi(x)|\le 2\lambda d(x,y)+2m$ for all $x,y\in X$.
\qed
\enddemo
REMARK 1.
In the proof of Theorem 4.1 in the case of $\sA$ we switched from a Lipschitz constant $\lambda$ to
$2\lambda$. We note that we can take $\alpha\lambda$ for any $\alpha>1$.

REMARK 2.
The space $\R_+$ would not be an $AE(\sA)$ if we take coarse maps in Roe's sense
as morphisms
in $\sA$. Then the map $\phi:A\to \R_+$ for $A=\{n^2\}\subset\R$ defined by
the formula $\phi(n^2)=2^n$, being a coarse map, does not have a coarse
extension, since no extension of $\phi$ can be asymptotically Lipschitz.

EXAMPLE 1. $\R^n$ is not an $AE(\sA)$. It is not an absolute extensor for itself.
To see that one can take $X=\R^n$
and $A=\cup S_{2^n}(0)$, the union of spheres of rapidly growing radii, and
define a map $f:A\to\R^n$ as the union of maps $f_n:S_{2^n}(0)\to S_n(0)$
of degree $n$. Then $f$ does not have an extension to a proper map.

EXAMPLE 2. A simpler proof that $\R^n$ is not an absolute extensor is 
the following: there is no proper retraction of $\R^{n+1}_+$ onto $\R^n$.
This proof is good for $\tilde\sA$ as well.

\proclaim{Lemma 4.2}
The space $\R^n_+$ is isomorphic in $\tilde\sA$ 
(and hence in $\sA$) to the $n$-th power $(\R_+)^n$.
\endproclaim
\demo{Proof}
There is the natural embedding of $(\R_+)^n$ in $\R^n$ generated by the
imbedding  $\R_+\subset\R$. We imbed $\R_+^n$ in $\R^n$ as the halfspace
bounded by the hyperplane $A$: $x_1+\dots+x_n=0$ which contains $(\R_+)^n$.
We define a morphism $T:(\R_+)^n\to\R_+^n$ by the following rule.
For any line $l$ parallel to the vector $\bar d=(1,\dots,1)$ the restriction $T_l$ of
$T$ onto $l\cap(\R_+)^n$ is a translation by vector $\bar w_l$ in the direction
$-\bar d$. The norm $\|\bar w_l\|$ we denote by $w_l$. If $\bar x=(x_1,\dots,x_n)$
is a point on $l$ then we define $w_l=\frac{1}{n}(\Sigma x_i-mn)$ where 
$m=\min{x_i}$.
We note that the definition of $w_l$ does not depend on the choice of 
$\bar x\in l$.
Note that $m=0$ for a point $\bar x\in\partial(\R_+)^n$. 
Since $\Sigma_j(x_j-\frac{1}{n}\Sigma_i x_i)=0$,
the boundary of $(\R_+)^n$ is taken to $A$ by $T$ and hence $T$ is surjective.
This implies that $T$ is bijective and $T^{-1}$ is the family of inverse
translations $\{-\bar w_l\}$. If $\bar x$ and $\bar y$ are two points lying on 
lines $l$ and $l'$, then $\|\bar w_l-\bar w_{l'}\|=|w_l-w_{l'}|=
\frac{1}{n}|\Sigma x_i-\Sigma y_i+m'n-mn|\le \frac{1}{n}\Sigma|x_i-y_i|+|m'-m|$.
Show that $|m'-m|\le\Sigma|x_i-y_i|$. Indeed, we may assume that $m<m'$.
Let $m$ is achieved on $i$. Then $|y_i-x_i|=y_i-m\ge m'-m$.
Then $\|\bar w_l-\bar w_{l'}\|\le \frac{n+1}{n}\Sigma|x_i-y_i|\le 
nd(\bar x,\bar y)$. Therefore $d(T(\bar x),T(\bar y))\le (n+1)d(\bar x,\bar y)$.
The same Lipschitz constant is good for $T^{-1}$.

Note that $\|T(x)\|=\|x+\bar w_l\|\ge \|x\|-w_l\ge (1-\frac{1}{\sqrt{2}})\|x\|$.
The last inequality is due to the fact that $w_l=\|pr_Ax\|\le\frac{1}{\sqrt{2}}\|x\|$.
This implies that $T$ has nonzero norm. By virtue of the inequality
$\|T^{-1}(x)\|\ge \|x\|$ the inverse map $T^{-1}$ also has a nonzero norm.
\qed
\enddemo
\proclaim{Theorem 4.3}
$\R^n_+\in AE(\sA)$ and $\R^n_+\in AE(\tilde\sA)$ for all $n$.
\endproclaim
\demo{Proof}
According to Lemma 4.2 it suffices to prove that for the $n$-th power of $\R_+$.
We fix a basis $\bar v_1,\dots,\bar v_n$ in $\R^n$ with 
$\bar v_1=(1,\epsilon,\dots,\epsilon),\ldots,\bar v_n=
(\epsilon,\dots,\epsilon,1)$ for some small $\epsilon$.
We define a projection $p_i:(\R_+)^n\to \R_+$ onto the $i$-th factor by the
formula $p_i(\bar z)=\bar z\cdot\bar v_i$. Geometrically,
$p_i$ is a projection onto $x_i$-axis, parallel to the plane $\alpha_i$,
where $\alpha_i$ is orthogonal to $\bar v_i$.
Note that each $p_i$ is a proper Lipschitz map. 

We define a linear isomorphism $p':\R^n\to\R^n$ with 
$p'(\bar e_i)=\bar v_i$ where $\{\bar e_i\}$ is the standard orthonormal basis 
in $\R^n$. It is easy to check that $p'(\bar z)=\Sigma p_i(\bar z)\bar e_i$. 
Denote by 
$p:(\R_+)^n\to p'((\R_+)^n)=V$  the restriction of $p'$ onto $(\R_+)^n$. 
The inverse map $p^{-1}:V\to(\R_+)^n$ can be extended to a Lipschitz map
$q:(\R_+)^n\to(\R_+)^n$. To demonstrate that first we note that $V$ is the union of rays in 
$\R^n$ emanated from the origin through a simplex $\sigma\subset\Delta^{n-1}$
lying in the standard simplex $\Delta^{n-1}=\{\bar x\mid\Sigma x_i=1, x_i\ge 0\}$.
Moreover the simplex $\sigma$ is obtained from $\Delta^{n-1}$ by a contraction
$c$ with the fixed point in the center of $\Delta^{n-1}$. Let 
$\gamma=c^{-1}:\sigma\to\Delta^{n-1}$ be the inverse map. We can extend $\gamma$
to $\xi:\Delta^{n-1}\to\Delta^{n-1}$ by taking the radial contraction
of the collar $\Delta^{n-1}\setminus\sigma$ to $\partial\Delta^{n-1}$.
The map $\xi$ can be linearly extended to a map $\beta:(\R_+)^n\to(\R_+)^n$.
It is clear that $\beta$ is Lipschitz. The restriction $\beta\mid_V$ is a linear 
map which takes $\bar v_i$ to $(1+(n-1)\epsilon)\bar e_i$. Hence,
$\frac{1}{1+(n-1)\epsilon}\beta\mid_V=p^{-1}$. Then 
$q=\frac{1}{1+(n-1)\epsilon}\beta$.

Let $\phi:A\to\R^n_+$
be a morphism in $\sA$ ($\tilde\sA$), where $A\subset Z$. Then
$p_i\circ\phi$ is a morphism in $\sA$ ($\tilde\sA$) for any $i$. 
By Theorem 4.1 there are
extensions $\bar\phi_i:Z\to\R_+$ with constants $\lambda_i$, $s_i$.
Consider the map $\psi=q\circ\bar\phi:Z\to(\R_+)^n$ where $\bar\phi=(\bar\phi_1,\dots,
\bar\phi_n)$. For $z\in A$ we have $\psi(z)=q(\bar\phi_1(z),\dots,\bar\phi_n(z))
=p^{-1}(p_1\circ\phi(z),\dots,p_n\circ\phi(z))=p^{-1}p\phi(z)=\phi(z)$.
Thus, $\psi$ is an extension of $\phi$. The map $\psi$ is proper 
(with nonzero norm), since
all $\bar\phi_i$ are of that kind. The map $\psi$ is asymptotically Lipschitz
with constants $\lambda m$ and $s$ where $m=\max\{\lambda_i\}$, $s=\max\{s_i\}$
and $\lambda$ is a Lipschitz constant for $q$.
\qed
\enddemo

\proclaim{Lemma 4.4}
Let $A$ be a closed subset of a proper metric space $(X,d)$ and let $g:A\to\R^n$
be an asymptotically Lipschitz map. Then there is a neighborhood $W\supset A$
with $\|g(a)\|\le \lambda(d(a,X\setminus W))+s$ for some numbers $\lambda$
and $s$
and for all $a\in A$ which admits an asymptotically Lipschitz 
extension $\bar g:W\to\R^n$ of $g$.
\endproclaim
\demo{Proof}
Let $\pi:\R^{n+1}_+\to\R^n$ be the orthogonal projection. By Theorem 4.3 there
is an asymptotically Lipschitz extension $g':X\to\R^{n+1}_+$ with constants
$\lambda'$ and $s'$. We define $\lambda=\sqrt{2}\lambda'$,
$s=\sqrt{2}s'$ and $W=(g')^{-1}(V)$ where $V$ is the region under
the graph of the function $\|x\|$ on $\R^n\subset\R^{n+1}_+$.
Then $\bar g=\pi\circ g'$.
Let $w\in X\setminus W$, then $\|g(a)-g'(w)\|\ge\frac{1}{\sqrt{2}}\|g(a)\|$.
Therefore, $\lambda'd(a,w)+s'\ge\frac{1}{\sqrt{2}}\|g(a)\|$.
Hence, $\lambda(d(a,w))+s\ge\|g(a)\|$ for all $w\in X\setminus W$.
Thus, $\lambda(d(a,X\setminus W))+s\ge\|g(a)\|$ for all $a\in A$.
\qed
\enddemo
REMARK 4.  Lemma 4.4 holds true if one replace $\R^n$ by $\R^n_+$.

We recall that a coarsely proper map $f$ is characterized by the
property that a preimage $f^{-1}(C)$ of any bounded set is
bounded.
\proclaim{Proposition 4.5}
For any coarsely proper (not necessarily continuous)
function $f:X\to \R_+$ of a proper metric space $(X,d)$
there is a proper asymptotically Lipschitz function $q:X\to\R_+$ with $q\le f$.
\endproclaim
\demo{Proof}
Denote $A_k=f^{-1}([0,k])$. Consider an increasing sequence
of concentric balls $B_{m_k}(x_0)$ of integral radii, $x_0\in X$ 
and $k\in\Bbb N$ such that
$A_k\subset B_{m_k}(x_0)$. Define $q(x)=k-2+\frac{d(x,x_0)-m_{k-1}}{m_k-m_{k-1}}$
for $x\in B_{m_k}(x_0)\setminus B_{m_{k-1}}(x_0)$ and $k>1$. We set
$q(B_{m_1}(x_0))=0$. Clearly $q$ is a proper continuous map. Since
$f(B_{m_k}(x_0)\setminus B_{m_{k-1}}(x_0))\ge k-1$ and
$q(B_{m_k}(x_0)\setminus B_{m_{k-1}}(x_0))\le k-1$, we have $q\le f$.
To check the asymptotical Lipschitz condition we take two points
$x\in B_{m_k}(x_0)\setminus B_{m_{k-1}}(x_0)$ and 
$y\in B_{m_l}(x_0)\setminus B_{m_{l-1}}(x_0)$ with $k\ge l$.
Then $|q(x)-q(y)|=|k-l+\frac{d(x,x_0)-m_{k-1}}{m_k-m_{k-1}}-
\frac{d(y,x_0)-m_{l-1}}{m_l-m_{l-1}}|\le |k-l|+2\le d(x,y)+3$.
\qed
\enddemo
\proclaim{Lemma 4.6}
Let $A$ be a closed subset of a proper metric space $(X,d)$,
let $W\supset A$ be a closed neighborhood and let $g:W\to\R_+$ be   proper asymptotically 
Lipschitz with $\lambda(d(a,X\setminus W))+s\ge g(a)$ for some constant $\lambda$
and $s$ and all $a\in A$. Let $f:X\to\R_+$ be a coarsely proper map with $g\le f$.
Then there exists a proper asymptotically Lipschitz map $\bar g:X\to\R_+$ with
$\bar g\le f$ and $\bar g\mid_A=g$.
\endproclaim
\demo{Proof}
Since $\R_+$ is AR, there is an asymptotically Lipschitz extension $g':X\to\R_+$
of $g$. Apply Proposition 4.5 to a proper map $f'=\min\{f,g'\}$ to obtain
an asymptotically Lipschitz map
 $q:X\to\R_+$ such that $q\le f$ on $X$ and $q\le g$ on $W$.
We define $\phi(x)=\frac{d(x,X\setminus W)}{d(x,A)+d(x,X\setminus W))}$.
Since $\phi(X\setminus A)=0$, the function $\bar g(x)=q(x)+\phi(x)(g(x)-q(x))$
is well-defined. Since $q\le \bar g$ and $\phi\ge 0$, we
have that $\bar g(x)$ tends to infinity when $x$
approaches infinity i.e. when $d(x,x_0)\to\infty$ for 
some (all) fixed point $x_0\in X$.
Note that $\bar g\mid_A=g$. To complete the proof we verify that $\bar g$ is
asymptotically Lipschitz. First, we note that

\

$|\phi(x)-\phi(y)|=\frac{1}{d(y,A)+d(y,X\setminus W)}
\frac{|d(y,A)d(x,X\setminus W)-d(x,A)d(y,X\setminus W)|}{d(x,A)+
d(x,X\setminus W)}\le$

\

$\frac{1}{d(y,A)+d(y,X\setminus W)}
(|d(y,A)-d(x,A)|\phi(x)+
|d(x,X\setminus W)-d(y,X\setminus W)|(1-\phi(x))$

\

$\le\frac{1}{d(y,A)+d(y,X\setminus W)}(d(x,y)
\phi(x)
+d(x,y)(1-\phi(x)))=
\frac{1}{d(y,A)+d(y,X\setminus W)}d(x,y)$. 

\

Then
$|\bar g(x)-\bar g(y)|\le |q(x)-q(y)+\phi(x)g(x)-\phi(x)q(x)-\phi(y)g(y)+\phi(y)q(y)|
\le|q(x)-q(y)|+\phi(x)|g(x)-g(y)|+g(y)|\phi(x)-\phi(y)|+q(y)|\phi(x)-\phi(y)|+
\phi(x)|q(x)-q(y)|\le 2|q(x)-q(y)|+|g(x)-g(y)|+2g(y)|\phi(x)-\phi(y)|$.
Here we use the facts that $q(y)\le g(y)$ and $\phi(x)\le 1$. Let $\lambda_1, s_1$
and $\lambda_2,s_2$ be constants from the definition of the asymptotical
Lipschitz property for $g$ and $q$. Then we can conclude that

\

$|\bar g(x)-\bar g(y)|\le 2\lambda_2d(x,y)+2s_2+\lambda_1d(x,y)+s_1+
2g(y)\frac{1}{d(y,A)+d(y,X\setminus W)}d(x,y)\le$

\

$(2\lambda_2+\lambda_1+
\frac{\lambda d(y,X\setminus W)+s}{d(y,A)+d(y,X\setminus W)})
d(x,y)+s_1+2s_2$. 

\

Let $m=\min\{d(y,A)+d(y,X\setminus W)\}$. Then
$\frac{\lambda d(y,X\setminus W)+s}{d(y,A)+d(y,X\setminus W)}\le(\lambda+s/m)$.

\

Therefore $|\bar g(x)-\bar g(y)|\le(\lambda_1+2\lambda_2+\lambda+s/m)d(x,y)+
s_1+2s_2$.\qed

\enddemo

\head \S5 ANE, Homotopy, Dimension \endhead

For any closed subset $A\subset X$ of a proper metric space $X$ we
define an {\it asymptotic neighborhood} $W$ of $A$ in $\sA$ as a subset of $X$ 
containing $A$ with
the property $\lim_{R\to\infty}d(A\setminus B_R(x_0),(X\setminus W)\setminus
B_R(x_0))=\infty$ for some (= any) point $x_0\in X$.

Similarly, in $\tilde\sA$, a set $W$ should have the property
$d(x, X\setminus W)\ge k\|x\|$ for some $k>0$ for all $x\in A$.

DEFINITION 1 (conventional). An object $Y\in\sA$ (or $\tilde\sA$) is called
an {\it absolute neighborhood extensor}, $Y\in ANE(\sA)$, if for any object
$X$ and any subobject $A\subset X$, for every morphism $f:A\to Y$ there is
an extension $\bar f:W\to Y$ of $f$ to a morphism of a closed 
asymptotic neighborhood.

This definition has a flaw. In the category of topological spaces,
metric compacta in particular, there is a fact connecting $AE$ and $ANE$.
Namely,  $X$ is ANE if and only if the cone over it $CX$ is AE. In category
$\sA$ this is not the case. Let $X$ be a parabola, lying in $\R^2_+$
with the induced metric. Since $X$ is a retract of an asymptotic neighborhood 
in $AE$-space $\R^2_+$, $X$ is a conventional ANE. It is not difficult to show
that the imbedding of $X$ in $CX$ cannot be extended over $\R^2_+$. If one prefers
to consider geodesic metric spaces, he could take a paraboloid with 
an inner metric.

Note that the implication $CX\in AE\ \ \ \Rightarrow \ \ \ X\in ANE$ always holds.

DEFINITION 2 (categorical). $X$ is $ANE_0$ if $X\times\R_+\in AE$.

Clearly $ANE_0(\sA)\subset ANE(\sA)$.
\proclaim{Proposition 5.1} $ANE_0(\tilde\sA)\subset ANE(\tilde\sA)$.
\endproclaim
\proclaim{Proposition 5.2} 
$\R^n\in ANE_0$ in both categories $\sA$ and $\tilde\sA$.
\endproclaim
\demo{Proof}
The result follows from Theorem 4.3.
\enddemo
To make the analogy with the topological category more visual we
denote the half plane $\R^2_+$ by $\Bbb I$.
\proclaim{Theorem 5.3 (HET)}
Let $Y\in ANE(\tilde\sA)$. Let $A$ be a closed subset of $X$ and
let $f:i_{-}(X)\cup A\tilde\times \Bbb I\to Y$ be a morphism. Then there
is an extension $\bar f:X\tilde\times \Bbb I\to Y$ to a morphism.
\endproclaim
\demo{Proof} 
Since $Y\in ANE$, there is an extension $g:W\to 
Y$ of $f$ to a neighborhood.  It implies that
$d(x,(X\tilde\times I)\setminus W)\ge \lambda\|x\|$.
Therefore one can construct a Lipschitz function $\phi:X\to\R_+$
extending a function $\psi$ on $A$, given by the formula $\psi(x)=d(x,x_0)+1$,
such that the region $D$ under the graph of this function lies in
$W$. 
Here we consider the realization of $X\tilde\times \Bbb I$
in the space $X\times\R_+$. Consider the map $h:X\times \Bbb I\to D$, defined by
taking each interval $\{x\}\times[0,d(x,x_0)+1]$ linearly to 
$\{x\}\times [0,\phi(x)]$. Let $\lambda'$ be a Lipschitz constant of
$\phi$. Show that $h$ is Lipschitz. Indeed,

\

$d_{X\times\R_+}(h(x,t),h(x',t'))=d_X(x,x')+|\frac{t}{\|x\|+1}\phi(x)-
\frac{t'}{\|x'\|+1}\phi(x')|\le $

\

$d_X(x,x')+
+\frac{t}{\|x\|+1}|\phi(x)-\phi(x')|+
\phi(x')|\frac{t}{\|x\|+1}-\frac{t'}{\|x'\|+1}|\le$ 

\

$d_X(x,x')+\lambda' d_X(x,x')+
\frac{\phi(x')}{\|x'\|+1}(t-t')+\frac{t}{(\|x\|+1)(\|x'\|+1)}|\|x'\|-\|x\| |\le$

\

$(\lambda'+1)d_X(x,x')+|t-t'|+d_X(x,x')\le (\lambda'+2)(d_X(x,x')+|t-t'|)=$

\

$(\lambda'+2)d_{X\times\R_+}((x,t),(x',t'))$.

\

Since $\|x\|+1\ge t$, we have $\|h(x,t)\|=\|x\|+\frac{t}{\|x\|+1}\phi(x)\ge
\frac{1}{2}(\|x\|+t)-\frac{1}{2}=\frac{1}{2}\|(x,t)\|-\frac{1}{2}$. Thus,
$h$ has nonzero norm. 

We define $\bar f= g\circ h$.\qed
\enddemo

DEFINITION. A {\it homotopy} between two morphisms $f,g:X\to Y$ is a morphism
$H:X\tilde\times \Bbb I\to Y$ such that $H\mid_{i_-(X)}=f$ and
$H\mid_{i_+(X)}=g$.

The notion of homotopy leads to the notion of homotopy equivalence.

EXAMPLE. As it was shown in [Roe2] $\R^n$ is homotopy equivalent to the
$n$-dimensional hyperbolic space $\H^n$ in $\sA$. They are not homotopy equivalent
in $\tilde\sA$, since $\R^n$ does not admit a Lipschitz degree one map
to $\H^n$ with nonzero norm.

DEFINITION. A metric space $X$ is said to be of {\it bounded geometry},
 $X\in BG$, if for every
$L$ there is a uniformly bounded cover $\sU$ of $X$ with the Lebesgue number 
$> L$ and of finite multiplicity.

Let $X$ be a metric space of bounded geometry, denote by $d(L)$ the minimal
multiplicity $m(L)$, which can be achieved in the above definition, minus one:
$d(L)=m(L)-1$.

DEFINITION 0 (Gromov. [G1]). The maximum $\max\{d(L)\mid L\in\R_+\}$, if exists, is
called the {\it asymptotic dimension} of $X$ and is denoted $as\dim X$.

This is a coarse analog of the Lebesgue covering dimension. In classical
topology there are several different definition of dimension which lead to
the same result (for nice class of spaces, say compact metric). We consider
here Ostrands' definition, Alexandroff-Urysohn's, and Alexandroff-Hurewicz'
(see $\S 1$).

DEFINITION 1 (Gromov [G1]). $as\dim X\le n$ if for any $L>0$ there are
$n+1$ $L$-disjoint families $\sU_i$ of uniformly bounded sets, $i=0,\dots, n$
such that the union $\cup \sU_i$ forms a cover of $X$.

This definition agrees with DEFINITION 0 [G1].

The {\it width} of a simplex $\Delta$ lying in a Banach space $(V,\|\ \|)$ is 
defined as the minimal distance from the barycenter of $\Delta$ to a face
$\sigma\subset\Delta$. 
An {\it asymptotic polyhedron} $P$ is a locally finite polyhedron supplied with an
intrinsic metric whose restriction on each simplex has a metric induced 
from some Banach space (the same for all simplices) such that widths of
simplices tend to infinity as one recede from a base point. The latter
means that for any $M>0$ there is a finite subcomplex $K\subset P$ such that
the width of any simplex from $P\setminus K$ is greater than $M$.

Let $\phi:X\to Y$ be a morphism and let $X,Y\subset Z$. The Alexandroff
norm of $\phi$ is the displacement function $\|\phi\|_A:X\to\R_+$ defined by 
the formula $\|\phi\|_A(x)=d_Z(x,\phi(x))$. 
Let $f:X\to\R_+$ be a proper function. We say that {\it the 
Alexandroff-Gromov norm}
of $\phi:X\to Y$ does not exceed $f$, $\|\phi\|_{AG}<f$,
 if there exists a metric space $Z$ and
isometric imbeddings $X,Y\subset Z$ such that $\|\phi\|_A(x)<f(x)$.
We also define {\it the Urysohn norm} of $\phi:X\to Y$ to be majorated by
a proper function $f:X\to\R_+$, $\|\phi\|_U<f$, if for any $R$ there is a compactum
$C\subset X$ such that $diam(\phi^{-1}(B_R(\phi(x))))<f(x)$ for all 
$x\in X\setminus C$. It is easy to see that $(2+\epsilon)\|\phi\|_{AG}>
\|\phi\|_U$ for any $\epsilon>0$.

DEFINITION 2. A metric space $X$ has an asymptotic dimension
$as\dim_*X\le n$ if for any proper function $f:X\to\R_+$
there is a short map $\phi:X\to K$ to an $n$-dimensional asymptotic polyhedron
such that $\|\phi\|_U<f$.

One can use here the Alexandroff-Gromov norm in this definition instead of the 
Urysohn norm. 

Let $\sU$ be a cover of a metric space $X$ and let $x\in X$. We denote

$L_{\sU}(x)=\sup\{d(x,X\setminus U)\mid U\in\sU\}$,

$mesh_{\sU}(x)=\sup\{diam(U)\mid U\in\sU, x\in U\}$,

$order_{\sU}(x)=m_{\sU}(x)=card\{U\in\sU\mid x\in U\}$.
\proclaim{Proposition  5.3}
$asdim_*X\le n\Leftrightarrow$ for a given proper map $f:X\to\R_+$ there is a uniformly
bounded cover $\sU$ of order $\le n+1$ with $\lim_{x\to\infty}L_{\sU}(x)=\infty$ and 
with $mesh_{\sU}(x)<f(x)$ for all $x\in X\setminus C$ for some compact set $C$.
\endproclaim
\demo{Proof}
1) Assume that $asdim_*X\le n$ and let $f:X\to\R_+$ be given. By the definition
there is a short map $\phi:X\to K$ to an asymptotic polyhedron with
$\|\phi\|_U< f$. We may assume that $\phi$ is onto. Let $C_R\subset X$ be
a compactum from the definition of the inequality $\|\phi\|_U< f$.
Consider a filtration $T_1\subset T_2\subset\dots$ of $K$ by
subcomplexes such that $\phi(C_i)\subset T_i$. One can define a subdivision $K'$
of $K$ such that $mesh(x)\le i/2$ for $x\in T_{i+1}$ and with width tending
to infinity. This is possible by taking the standard cubification, regular 
subdivision of cubes and then barycentric triangulation. We define 
$\sU=\{\phi^{-1}(OSt(v,K'))\mid v\in(K')^{(0)}\}$ where $OSt(v,K')$ means
the open star of vertex $v$ in $K'$. 
Assume that $v\in T_{i+1}\setminus T_i$.
Note that $diam\phi^{-1}(OSt(v,K'))\le
diam\phi^{-1}(B_i(y))<f(x)$ for all $y\in OSt(v,K')$ and $x\in\phi^{-1}(y)$. 
Thus, $mesh_{\sU}(x)<f(x)$. The Lebesgue number $L_{\sU}(x)$ can be estimated
from below as minimum of width of $n$-simplices $\Delta$ with $\Delta\cap\Delta_x
\ne\emptyset$ where $\Delta_x$ contains $\phi(x)$. Thus 
$\lim_{x\to\infty}L_{\sU}(x)=\infty$.

2) Let $f$ be given. We may assume that $f(x)\le\|x\|$ where $\|x\|$ is the distance
to the base point $x_0$ in $X$. Let $g(t)=\inf\{f(x)\mid x\in X\setminus B_t(x_0)\}$
and let $\bar f(x)=\frac{1}{4}g(\|x\|-f(x))$. Not that $\bar f$ tends to infinity
as $x$ approaches infinity.
Consider a uniformly bounded cover $\sU$ of order $\le n+1$ with
$\lim_{x\to\infty}L_{\sU}(x)=\infty$ and with $mesh_{\sU}(x)<\bar f(x)$ 
for $x\in X\setminus C$ for some compact $C$. Let $\phi:X\to N(\sU)$ be projection
to the nerve of $\sU$. We can take a piecewise Euclidean metric on $N(\sU)$
which turns $N(\sU)$ into an asymptotic polyhedron and which makes the map
$\phi$ short. We have to show that for any $R$ for $x$ close enough to infinity
the inequality $diam(\phi^{-1}(B_R(\phi(x))))<f(x)$ holds. First we note
that for $x$ close enough to infinity the ball $B_R(\phi(x))$ lies in the star 
$St(\{V\},N(\sU))$, where $x\in V\in\sU$. Hence 
$\phi^{-1}(B_R(\phi(x)))\subset\cup_{U\cap V\ne\emptyset,U\in\sU}U=V^*$.
Let $y\in U$, $U\cap V\ne\emptyset$. Then $d(x,y)\le mesh_{\sU}(x)+mesh_{\sU}(y)
< \bar f(x)+\bar f(y)=\frac{1}{4}(g(\|x\|-f(x))+g(\|y\|-f(y)))\le
\frac{1}{4}(f(x)+\inf\{f(z)\mid z\in X\setminus B_{\|y\|-f(y)}(x_0)\})$.
The triangle inequality implies that $\|x\|\ge\|y\|-f(y)$. Hence
$x\in B_{\|y\|-f(y)}(x_0)$ and hence, 
$\inf\{f(z)\mid z\in X\setminus B_{\|y\|-f(y)}(x_0)\}\le f(x)$. We conclude
that $d(x,y)<\frac{1}{2}f(x)$. Therefore $diam(V^*) <f(x)$. \qed
\enddemo

The method of proof of Theorem 1.1 of [D-K-U] allows to proof the following.
\proclaim{Proposition 5.4} $as\dim_*X\le as\dim X$ for all $X$.
\endproclaim
The following definition of a macroscopic dimension is based on the idea
of extending of maps to spheres.

DEFINITION 3. A coarse dimension does not exceed $n$, $\dim^cX\le n$, if for any closed subspace
$A\subset X$, any morphism $f:A\to\R^{n+1}$ in $\sA$, there is an extension
$\bar f:X\to \R^{n+1}$.

Here $\R^{n+1}$ serves as an analog of $n$-sphere in $\sA$. There reason to
believe that this is a proper analog is based on the following. An analog of
a point (minimal AE-space) in $\sA$ is a halfline $\R_+$. Then an analog of
0-dimensional sphere $S^0$ might be $\R=\R_+\cup\R_+$. The suspension in
this case is a multiplication by $\R$.

Latter we will prove the inequality $\dim^cX\le as\dim_*X$.

DEFINITION. A metric space $X$ has a {\it slow dimension growth} if
$\lim_{L\to\infty}\frac{d(L)}{L}=0$.

\head \S6 Approximation by Polyhedra, Cohomology and Dimension\endhead

A proper metric space $X$ of bounded geometry admits so called {\it anti-\v{C}ech
approximation} by polyhedra. This is by the definition a direct sequence
$\{K_i,g^i_{i+1}:K_i\to K_{i+1}\}$ with simplicial projections $g^i_{i+1}$
together with short maps $f_i:X\to K_i$ such that $f_{i+1}$ is homotopic to
$g^i_{i+1}\circ f_i$. All simplices in $K_i$ are isomorphic
to the standard simplex of size $L_i$, the metric on $K_i$ is a geodesic metric, 
induced by this property, and $\lim L_i=\infty$. Moreover, one may assume that
there are short projections $f^i_{i+1}:K_i\to K_{i+1}$ such that 
$f_{i+1}=f^i_{i+1}\circ f_i$ and $g^i_{i+1}$ are simplicial approximations of 
$f^i_{i+1}$. When $X$ is uniformly contractible, maps $f_i$ admit left proper homotopy
inverse $p_i:K_i\to X$ such that $f_{i+1}^i$ and $f_{i+1}\circ p_i$ are 
properly homotopic (see [Roe1], [H-R] for more details).

Here we consider simplicial complexes of the kind that participate in
the definition of an anti-\v{C}ech approximation, i.e. simplicial complexes
with a geodesic metric such that all simplices are isomorphic to the standard
simplex of size $L$. We refer to $L$ as the mesh of $K$, $mesh(K)$ and we call 
such simplicial complexes as {\it uniform polyhedra}.

Following Roe [Roe1], [H-R], we define anti-\v Cech homology (coarse homology)
of $X$ as $\hat H_*(X)=\lim_{\to}\{H_*(K_i)\}$. One can take homologies with
infinite locally finite chains, as the result he will get Roe's exotic homology
$\hat H^{lf}_*(X)=HX_*(X)$ as in [H-R].  This definition is dual to that
for classical \v Cech homology and cohomology for compact metric spaces.
We recall that in the case of compacta the \v Cech homology generally does not 
behave nicely. It is not exact. There is an exact homology theory called
the Steenrod homology. If a compactum $Y$ is the limit of an inverse sequence
of polyhedra $\{L_i,q^{i+1}_i\}$, the $k$-th Steenrod homology can be 
defined by the formula $H^s_k(Y)=H^{lf}_{k+1}(T)$ where $T$ is the telescope
generated by the inverse sequence. The Steenrod theory has a flaw on its own:
it is not continuous. So both theories are needed. Note when the
coefficient group is a field, then these two theories coincide.

In the large scale topology we also form a telescope $T$ out of the direct
sequence $\{K_i,g^i_{i+1}\}$ giving an Anti-\v Cech approximation to a
metric space $X$. Dual to the Steenrod homology we define the Roe cohomology
as $H^k_R(X)=H^k(T)$. To define Roe's cohomology with compact support we consider
the telescope $ T^{\alpha}$ of one point compactifications $\alpha K_i$ of $K_i$
generated by maps $g^i_{i+1}$. Then $H^k_{R,c}(X)=H^k(T^{\alpha})$. In Roe's
notations $H^k_{R,c}(X)=HX^k(X)$. There is a short exact sequence [Roe1]

$0\to\lim^1_{\leftarrow}H^{k-1}_c(K_i)\to H^k_{R,c}(X)\to \hat H^k_c(X)\to 0$.

A closed subset $A\subset X$ is {\it closed in category} $\sA$ if $X\setminus A$ is
an asymptotic neighborhood of some set in $X$. For every such $A$ the pair
$(X,A)$ has an anti-\v Cech approximation $\{(K_i,L_i),g^i_{i+1}\}$.
Similarly one can define homologies and cohomologies of a pair
$\hat H^{lf}_i(X,A)$, $\hat H^i_c(X,A)$ and $\hat H^i_{R,c}(X,A)$.
\proclaim{Proposition 6.1}
For any categorically closed subset $A\subset X$ there are exact sequences 
of a pair $(X,A)$

1) $\dots\to\hat H^{lf}_i(A)\to\hat H^{lf}_i(X)\to
\hat H^{lf}_i(X,A)\to\hat H^{lf}_{i-1}(A)\to\dots$ ,

\

2) $\dots\leftarrow H^i_{R,c}(A)\leftarrow H^i_{R,c}(X)
\leftarrow H^i_{R,c}(X,A)\leftarrow H^{i-1}_{R,c}(A)\leftarrow\dots$ .
\endproclaim
\demo{Proof}
Note that $\hat H^{lf}_k(X)=\lim_{\to}H^{lf}_k(K_i)=
\lim_{\to}H^s_k(\alpha K_i)=
\lim_{\to}H^s_k( T^{\alpha}_i)=H^s_k( T^{\alpha})$,  
where $T^{\alpha}_i$ is a finite
part of the telescope $ T^{\alpha}$ up to $\alpha K_i$. Then both exact sequences follow
from Steenrod homology and cohomology exact sequences 
of the pair $(T^{\alpha}_X,T^{\alpha}_A)$.
\enddemo
The inclusion of $\alpha X$ into a telescope $\alpha T_X$ generates homomorphisms
$c_{\ast}:H^{lf}_{\ast}(X)\to\hat H^{lf}_{\ast}(X)$ and $c^*: H^*_{R,c}(X)\to H^*_c(X)$.
\proclaim{Theorem 6.2[H-R],[Roe1]}
If $X$ is uniformly contractible, then $c_{\ast}$ and $c^*$ are isomorphisms.
\endproclaim
In the relative case we have the following
\proclaim{Theorem 6.3}
If $X$ is uniformly contractible, then for any categorically closed subset
$A\subset X$ there is the equality

\

1) $\hat H^{lf}_{\ast}(X,A)=\lim_{@>>k>}H^{lf}_{\ast}(X,N_k(A))$ 

\

and the exact sequence

\

2) $0\to\lim^1_{@<<k<}H^{\ast-1}_c(X,N_k(A))\to H^*_{R,c}(X,A)\to\lim_{@<<k<}H^*_c(X,N_k(A))\to 0$,

\

where $N_k(A)$ is the closed $k$-neighborhood of $A$.
\endproclaim
\demo{Proof}
Let $j^k_{k+1}:N_k(A)\to N_{k+1}(A)$ denote the inclusion.
We consider the telescope $S$ formed by these inclusions 
and let $S^{\alpha}$ be a corresponding telescope of the
one point compactifications $\alpha N_k(A)$. There is a natural
inclusion $S^{\alpha}\subset \alpha X\times\R_+$. We show that the pair
$(\alpha X\times\R_+,S^{\alpha})$ is proper homotopy equivalent to the
pair $(T^{\alpha}_{ X},T^{\alpha}_A)$. We may assume that
an anti-\v Cech approximation is chosen such that $f_i(N_i(A))\subset L_i$.
We define a map $f:X\times\R_+\to T_X$ as the union 
$\cup f_i$ on the set
$X\times\N$ and extend it over the rest of $X\times\R_+$ by means of homotopies
between $f^i_{i+1}$ and $g^i_{i+1}$. For every $i$ there is $k(i)$ such that
$p_i(L_i)\subset N_{k(i)}(A)$ and $p_i\circ f_i\mid_{N_{k(i-1)}(A)}$ is homotopic to
the identity $id_{N_{k(i-1)}(A)}$ in $N_{k(i)}(A)$. Then we define $p:T_X\to X\times\R_+$
as an extension of the union of maps $p_i:K_i\to X\times{k(i)}$.
Clearly $f$ and $p$ define a proper homotopy equivalence between $X\times\R_+$
and $T_X$. Easy verification shows that they also define a proper homotopy
equivalence of pairs $(X\times\R_+,S)$ and $(T_X,T_A)$. Then they induce a
proper homotopy equivalence of pairs $(T^{\alpha}_{ X},T^{\alpha}_A)$ and
$(\alpha X\times\R_+,S^{\alpha})$. Then
$\hat H^{lf}_i(X,A)=H^s_i(T^{\alpha}_X,T^{\alpha}_A)=H^s_i(\alpha X\times\R_+,S^{\alpha})=
\lim_{\to}H_i^{lf}(X,N_k(A))$ and
$H^*_{R,c}(X,A)=H^*(T^{\alpha}_X,T^{\alpha}_A)=H^*(\alpha X\times\R_+,S^{\alpha})$.
Then by Milnor's formula

$0\to\lim^1_{@<<k<}H^{\ast-1}_c(X,N_k(A))\to H^*_{R,c}(X,A)\to\lim_{@<<k<}H^*_c(X,N_k(A))\to 0$,

\qed
\enddemo
The above definitions of anti-\v Cech homology and
anti-\v Cech and Roe cohomologies works well for any
generalized homology and cohomology theory:
$\hat M^{lf}_{\ast}(X,A)=M^s_{\ast}(T^{\alpha}_X,T^{\alpha}_A)$,
$M^*_{R,c}(X,A)=M^{\ast}(T^{\alpha}_X,T^{\alpha}_A)$,
and $\hat M^*(X,A)=\lim_{\leftarrow}M^*(K_i,L_i)$.
The statements 6.1-6.3 also hold true for generalized
homology (cohomology) case.

DEFINITION. The asymptotical cohomological dimension of a proper metric space
$X\in BG$ with respect to coefficient group $G$ is 

$as\dim_GX=\sup\{n\mid H^n_{R,c}(X,A;G)\ne 0, \text{for some A}\}$.

\proclaim{Theorem 6.4}
For uniformly contractible space $X$ the inequalities 
$as\dim_GX\le\dim_GX\le\dim X$ hold,
where $\dim_G$ is the ordinary cohomological dimension
and $\dim$ is the classical covering dimension.
\endproclaim
\demo{Proof}
Follows from Theorem 6.3.
\enddemo
Gromov's Lemma 7.1 implies the following
\proclaim{Theorem 6.5}
A space $X$ has $as\dim X\le n$ if and only if
$X$ admits an anti-\v Cech approximation by
$n$-dimensional simplicial complexes.
\endproclaim

\head \S7 Higson corona \endhead

Let $f:X\to\R$ be a continuous function on a metric space $X$ and let
$B_R(x)$ be a ball centered at $x$ of radius $R$. The $R$-variation
$Var_Rf(x)$ of $f$ is the number $\sup_{y\in B_R(x)}|f(x)-f(y)|$.
Let $C_h(X)$ be a family of all bounded continuous functions on $X$
with the property $\lim_{x\to\infty}Var_Rf(x)=0$ for every $R$.
The {\it Higson compactification} of $X$ is the closure of 

$X=(\Pi_{f\in C_h(X)})f(X)\subset\R^{C_h(X)}$ in $\R^{C_h(X)}$.
The remainder $\nu X=\bar X\setminus X$ is called the {\it Higson corona}.
We note that the Higson corona is a covariant functor
$\nu:\sC\to Comp$ form the coarse category to the category of compact Hausdorff
spaces [Roe1].

A piecewise Euclidean (PE) complex of {\it mesh} $D$ is a simplicial complex $K$
whose simplices are isomorphic to the standard simplex of diameter $D$ and
the metric on $K$ is the geodesic metric induced by the metrics on simplices.
A map $f:X\to Y$ between metric spaces is called a {\it short map} if
$d_Y(f(x),f(y))\le d_X(x,y)$ for all $x,y\in X$. A map $f:X\to Y$ is called
{\it uniformly cobounded} if for any $R$ there is a constant $C$ such that the diameter
of the preimage $f^{-1}(B_R(y))$ does not exceed $C$ for any point $y\in Y$.
We recall that $asdim X\le n$ means for arbitrarily large $L>0$ there is
a uniformly bounded open cover of $X$ of multiplicity $\le n+1$
with the Lebesgue number $> L$.
\proclaim{Lemma 7.1 [G1]} For a proper metric space $X$ the following two
conditions are equivalent:
\roster
\item{} $asdim X\le n$;
\item{} For any $D>0$ there is a uniformly cobounded short proper map
$f:X\to K$ to a PE complex of mesh $D$ and dimension $n$.
\endroster
\endproclaim

\proclaim{Theorem 7.2}
If $asdim X<\infty$, then $asdim X=\dim\nu X$.
\endproclaim
\demo{Proof}
By [D-K-U], we have $asdim X\ge\dim\nu X$.

Let $asdim X=m$, we show that $\dim\nu X\ge m$. By Lemma 7.1 there is a sequence
of short maps $\phi_i:X\to K_i$ with $\dim K_i=m$ and $D_i=mesh(K_i)\to\infty$.
We define the Lipschitz constant of a family of maps $S=\{f:X\to Y\}$ as
$L(S)=\inf_{f\in S}\{\lambda\mid d_Y(f(x),f(x'))\le\lambda d_X(x,x')\}$
. If this number is not defined we set $L(S)=\infty$.
Let $g:X\to B^n$ be an inessential map, i.e. a map which admits a sweeping
$q:X\to\partial B^n$ with 
$q\mid_{g^{-1}(\partial B^n)}=g\mid_{g^{-1}(\partial B^n)}$.
Let $S(g)$ denote the set of all continuous sweepings of $g$. Fix a base point
$x_0\in X$ and define $\lambda^r_i=\sup\{L(S(\phi_i\mid_{\phi_i^{-1}(\Delta^m)}))
\mid \Delta^m\subset K_i\setminus IntB_r(\phi_i(x_0))\}$. Let $\lambda_i=
\overline{\lim_{r\to\infty}}\lambda_i^r$.

Show that that the sequence $\frac{D_i}{\lambda_i}$ is bounded. Assume not,
then there is a subsequence with $\frac{D_{i_k}}{\lambda_{i_k}}>k$.
This implies the inequality $\frac{D_{i_k}}{\lambda_{i_k}^{r_k}}>k$ for some
$r_k$. Hence $\lambda^{r_k}_{i_k}<\infty$ and therefore the maps
$\phi_{i_k}\mid_{\phi_{i_k}^{-1}(\Delta^m)}:\phi_{i_k}^{-1}(\Delta^m)\to\Delta^m$
are inessential for all $m$-simplices $\Delta^m\subset K_{i_k}\setminus 
IntB_{r_k}(\phi_{i_k}(x_0))$.
For every such simplex $\Delta^m\subset K_{i_k}$ we consider a sweeping
$\psi_{\Delta^m}:\phi_{i_k}^{-1}(\Delta^m)\to\partial\Delta^m$ with the
Lipschitz constant equal to $L(S(\phi_{i_k}\mid_{\phi_{i_k}^{-1}(\Delta^m)}))$.
The union of such sweepings defines a uniformly cobounded map
$\psi_{k}:X\to K_{i_k}^{(m-1)}\cup L_k$ where $L_k$ is the star neighborhood of
$B_{r_k}(\phi_{i_k}(x_0))$. Let $p:K_{i_k}^{(m-1)}\cup L_k\to M_k$ be the
simplicial map induced by collapsing $L_k$ to a point. Then $\dim M_k=m-1$.
Moreover the map $p$ is short if we consider PE structure on $M_k$ of mesh
$D_{i_k}$. If we rescale the metric on $M_k$ multiplying it by 
$\frac{1}{\lambda_{i_k}^{r_k}}$, then the composition $p\circ\psi_k$ will be
a short map as well. Note that $mesh(M_k)=\frac{D_{i_k}}{\lambda_{i_k}^{r_k}}>k$.
Lemma 7.1 implies that $asdim X\le m-1$, which contradicts the assumption.

Let $\frac{D_i}{\lambda_i}< b$ for all $i$. By induction we define a sequence
of $m$-simplices $\Delta_i^m\subset K_i$ such that the sets 
$A_i=\phi^{-1}_i(\Delta^m_i)$ are disjoint and 
$L(S(\phi_i\mid_{\phi_i^{-1}(\Delta^m)}))\ge
\frac{D_i}{4b}$. If $\Delta_1,\dots,\Delta_i$ are already defined, we take
$r_0$ such that $\cup_{l=1}^{l=i}A_l\subset\phi^{-1}_{i+1}(B_{r_0}
(\phi_{i+1}(x_0)))$. Then there is $r>r_0$ with $\lambda_{i+1}^r>\frac{D_{i+1}}{2b}$.
Take $\Delta^m_{i+1}\subset K_{i+1}\setminus IntB_r(\phi_{i+1}(x_0))$ with
$L(S(\phi_{i+1}\mid_{\phi_{i+1}^{-1}(\Delta^m)}))\ge\frac{D_{i+1}}{4b}$.

Let $p_i:\Delta^m_i\to\Delta^m$ be a linear map to the standard unit simplex.
Thus, $p_i$ is $D_i$-contraction. Let $A=\cup_{i=1}^{\infty}A_i$. Consider
a map $f:A\to\Delta^m$, defined by the formula 
$f=\cup_{i=1}^{\infty}p_i\circ\phi_i\mid_{A_i}$. Since $D_i$ tends to infinity, the
variation $Var_Rf(x)$ tends to zero for any given $R$. Hence there is an
extension $\bar f:\bar A=A\cup\nu A\to\Delta^m$. We show that the restriction
$\bar f\mid_{\nu A}$ is an essential map. Assume that it is inessential,
let $f':\nu A\to\partial\Delta^m$ be its sweeping. Denote 
$\partial A_i=\phi^{-1}_i(\partial\Delta^m)$. Consider the map
$g=f'\cup f\mid_{\cup_{i=1}^{\infty}\partial A_i}:\cup\partial A_i\cup\nu A\to
\partial\Delta^m$. Since $\cup\partial A_i\cup\nu A$ is a closed subset of
$\bar A$, there is an extension $\bar g:N\to\partial\Delta^m$ of $g$ over a neighborhood.
Then $A_i\subset N$ for large enough $i$. Let $g_i=\bar g\mid_{A_i}$.
Since the map $\bar g$ is extendible over $\nu A$, we have
$\lim_{i\to\infty}L(\{g_i\})=0$. On the other hand, $L(\{g_i\})\ge
L(S(p_i\circ\phi_i|_{\phi_i^{-1}(\Delta^m)}))\ge\frac{1}{D_i}\frac{D_i}{4b}=
\frac{1}{4b}$. The contradiction completes the proof.\qed
\enddemo
\proclaim{Lemma 7.3}
If $as\dim_*X<\infty$, then $as\dim_*X=\dim\nu X$
\endproclaim
\demo{Proof} 
First we give one more reformulation of the definition $asdim_*$. Namely,
$asdim_*X\le m\ \ \Leftrightarrow$ for any proper function $f:X\to\R_+$ there
is a map $\phi:X\to K$ to a locally finite piecewise euclidean simplicial complex of mesh one and of dimension $\dim K=m$ such that
$\lim_{\Delta}L(\{\phi\mid_{\phi^{-1}(\Delta)}\}=0$ and
$diam\phi^{-1}(\Delta)<\min_{x\in\phi^{-1}(\Delta)}f(x)$ for all simplices $\Delta\subset K\setminus C$ for some compact set $C$. We use notations taken from the proof of Theorem 7.2, i.e. $L(\{\phi\mid_{\phi^{-1}(\Delta)}\})$ is
the minimal Lipschitz number for the map $\phi:\phi^{-1}(\Delta)\to\Delta$.
Also we will write $\phi\prec f$ for the above condition on fibers of $\phi$.
It is not difficult to derive this equivalence from Proposition 5.3.

Let $as\dim_*X=m$, show that $\dim\nu X\ge m$.

Claim: There exist a proper map $\phi:X\to K$ onto PE complex of $\dim K=m$ with
mesh one and a sequence of $m$-simplices $\Delta^m_k$ with
$L(S\{\phi\mid_{\phi^{-1}(\Delta^m_k)}\})>a$ for some numbver $a>0$. Indeed,
if we assume the contrary then for every proper monoton function
$f:\R_+\to\R_+$, $f\le\frac{1}{2}t$, there is a map $\phi:X\to K$ to a complex
$K$ as above with $\phi\prec\frac{1}{2}f(\frac{1}{2}\|\ \|)$. Then for every simplex $\Delta^m\subset K$ we take a sweeping $\psi_{\Delta^m}$ with the lowest
possible Lipschitz number. Then the union $\psi=\cup_{\Delta\in K}\psi_{\Delta^m}:X\to K^{m-1}$ will be a map to a PE complex of dimension $m-1$
and of mesh one with $\lim_{\sigma\subset K^{(m-1)}}L(\{\psi\mid_{\psi^{-1}(\sigma)}\})=0$. Note that $\psi^{-1}(\sigma)\subset\cup_{\sigma\subset\Delta}\phi^{-1}(\Delta)$.
Then $diam\psi^{-1}(\sigma)\le 2\max_{\sigma\subset\Delta}diam\phi^{-1}(\Delta)\le 2\max_{\sigma\subset\Delta}\min\{\frac{1}{2}f(\frac{1}{2}\|x\|)\mid x\in\phi^{-1}(\Delta)\}\le\max\{f(\frac{1}{2}\|x\|)\mid x\in\cup_{\sigma\subset\Delta}\Delta\}=M$. Assume that $M=f(\frac{1}{2}\|x_1\|)$
for $x_1\in\cup_{\sigma\subset\Delta}\Delta$. Let $m=\min_{x\in\psi^{-1}(\sigma)}f(\|x\|)=f(\|x_2\|)$. Since
$diam\psi^{-1}(\sigma)\le M$, it follows that $\|x-2\|\ge\|x_1\|-M$.
Since $f(t)\le\frac{1}{2}t$, we have $\|x_2\|\ge\frac{1}{2}\|x_1\|$.
Since $f$ is monotone, we obtain $m\ge M$. Thus,
$diam\psi^{-1}(\sigma)<\min_{x\in\phi^{-1}(\sigma)}f(\|x\|)$. Since $f$ is
arbitrary enough, we have that $asdim_*X\le m-1$. Contradiction.

Assume that $\phi:X\to K$ and $\{\Delta^m_k\}$ as above. Denote
$A_k+\phi^{-1}(\Delta^m_k)$, and $A=\cup A_k$. Let $p_k:\Delta^m_k\to\Delta^m$
be the identity map onto the standard simplex. The union map $g=\cup_kp_k\circ\phi|mid_{A_k}:A\to\Delta^m$ can be extended to a map
$\bar g\bar A=A\cup\nu A\to\Delta^m$. Show that $\bar g\mid_{\nu A}$ is
essential. Assume the contrary: there is sweeping $g':\nu A\to\partial\Delta^m$.
Then $g'$ is extendible to a sweeping $\tilde g:\cup_{k\ge l}A_k\cup\nu A\to
\partial\Delta^m$. In this case $L(S\{\phi\mid_{A_k}\})\to 0$. That contradicts
with the inequality $L(S\{\phi\mid_{A_k}\})>a>0$.\qed
 \enddemo
Lemma 7.3 implies the following
\proclaim{Theorem 7.4}
Either $as\dim_*X=\dim\nu X$ or $as\dim_*X=as\dim X$.
\endproclaim
\proclaim{Proposition 7.5}
Let $f_n:X\to\R_+$ be a sequence of coarsely proper functions on a proper metric
space with $f_n\mid_W\ge g$ for some function $g:W\to\R$, given on a closed subset 
$W\subset X$, and for all $n$. 
Then there is a sequence of bounded 
subsets $A_n\subset X$ and a coarsely proper 
function $f:X\to\R_+$ with $f\mid_W\ge g$ such that $f\mid_{A_n}\le n$ and 
$f\mid_{X\setminus A_n}\le f_n$.
\endproclaim
\demo{Proof}
Let $\bar g:X\to\R_+$ be an extension of $g$ with the property $\bar g\le f_n$ for all
$n$.
Define $A_n=\cup_{i=1}^{i=n+1}f^{-1}_i([0,n])$ 
and $f(x)=\min\{\max\{n-1,\bar g(x)\},f_n(x)\}$ for
$x\in A_n\setminus A_{n-1}$. Note that $f_n(A_n\setminus A_{n-1})>n-1$ by virtue of
the definition of $A_{n-1}$. Hence $f(A_n\setminus A_{n-1})\ge n-1$. This implies
that $f(x)\to\infty$ as $x\to\infty$. Hence $f$ is coarsely proper.
If $x\in A_n$, then $x\in A_k\setminus A_{k-1}$ for some $k\le n$. 
Then $f(x)=\min\{\max\{k-1,\bar g(x)\},f_k(x)\}\le \max\{k-1,\bar g(x)\}$. 
Since $x\in A_k$, there is $i\le k+1$ such that $x\in f^{-1}_i([0,k])$.
Hence $\bar g(x)\le f_i(x)\le k$. Thus, $f(x)\le k\le n$.
If $x\in X\setminus A_n$, then $x\in A_m\setminus A_{m-1}$ for some
$m>n$. Since $x$ does not belong to $A_{m-1}$, we have $f_n(x)>m-1$.
Since $f_n(x)\ge\bar g(x)$, it follows that $f_n(x)\ge\max\{m-1,\bar g(x)\}$.
Therefore $f_n(x)\ge f(x)=\min\{\max\{m-1,\bar g(x)\},f_m(x)\}$.
Clearly, $f(x)\ge g(x)$.\qed
\enddemo
\proclaim{Theorem 7.6}
Let $X$ be a proper metric space. Then the following two conditions
are equivalent

1) $\dim^cX\le n$;

2) $\dim\nu X\le n$.
\endproclaim
\demo{Proof}
Assume that $\dim^cX\le n$. Let $\phi:C\to S^n$ be a map of a closed subset
of $\nu X$ to the unit $n$-sphere. There is an extension $\phi':V\to S^n$ of $\phi$
over a closed neighborhood in $\bar X=X\cup\nu X$. Then $Var_R\phi'(x)\to 0$
as $x\to\infty$ for any fixed $R$. We define proper functions
$f_n(x)=\frac{1}{Var_n\phi'(x)}$, $n\in\Bbb N$. Apply Proposition 7.5 to the
sequence $f_n$ with $g$ 
equal the constant zero function to obtain
a coarsely proper function $f:X\to\R_+$ and a filtration $X=\cup A_n$
such that $f(A_n)\le n$ and $f\mid_{X\setminus A_n}\le f_n$.

By Proposition 4.5 there is
an asymptotically Lipschitz function $q:X\to\R_+$ with $q\le f$
with Lipschitz constants $\lambda$ and $s$.
We define a map $g:X\cap V\to\R^{n+1}$ by the formula $g(x)=(q(x),\phi'(x))$
in the polar coordinates. We check that $g$ is asymptotically Lipschitz.
We denote by $\alpha(z_1,z_2)$ the angle between $z_1$ and $z_2$.
By the Law of Cosines we have

\

$\|g(x)-g(y)\|^2=q(x)^2+q(y)^2-2q(x)q(y)\cos(\alpha(\phi'(x),\phi'(y)))\le
(q(x)-q(y))^2+q(x)q(y)(\alpha(\phi'(x),\phi'(y)))^2$. 

\

Let $n-1\le d(x,y)\le n$.
Since $q(x)\le f_n(x)$ and $q(y)\le f_n(y)$ on $X\setminus A_n$, we have

\

$q(x)q(y)(\alpha(\phi'(x),\phi'(y)))^2\le f_n(x)f_n(y)(\alpha(\phi'(x),\phi'(y)))^2=
\frac{\alpha(\phi'(x),\phi'(y))}{Var_n\phi'(x)}
\frac{\alpha(\phi'(x),\phi'(y))}{Var_n\phi'(y)}$

$\le 1$. 

\

 For $x\in A_n$ we have $n\ge f(x)\ge q(x)$.
Hence, 

\

$q(x)q(y)(\alpha(\phi'(x),\phi'(y)))^2\le n(n+\lambda d(x,y)+s)\pi^2$. 

\

Therefore
$\|g(x)-g(y)\|^2\le(\lambda d(x,y)+s)^2+\pi^2((\lambda+1)d(x,y)+s+1)^2\le$

\

$(\lambda d(x,y)+s+\pi((\lambda+1)d(x,y)+s+1))^2$. 

\

Hence
$\|g(x)-g(y)\|\le (\lambda+\pi(\lambda+1))d(x,y)+(s+\pi(s+1))$.

\

Similarly if $y\in A_n$, then $\|g(x)-g(y)\|\le 
(\lambda+\pi(\lambda+1))d(x,y)+(s+\pi(s+1))$.

\

By the assumption there is an asymptotically Lipschitz extension
$\bar g:X\to\R^{n+1}$ of $g$. 
Hence $\bar g$ can be (uniquely) extended over Higson compactifications to
$\tilde g:\bar X\to\overline{\R^{n+1}}$
Consider the restriction
$g'=\tilde g\mid_{\bar X\setminus\tilde g^{-1}(0)}:\bar X\setminus\tilde g^{-1}(0)\to\overline\R^{n+1}\setminus\{0\}
$. 
Let
$\eta:\R^{n+1}\setminus\{0\}\to S^n$ be the radial projection 
and let $\bar\eta$ denote the extension to the Higson corona
$\bar\eta:\overline{\R^{n+1}}\setminus\{0\}\to S^n$.
Note that $\eta\circ g'$ restricted to $V\cap X$ coincides with $\phi'$.
Hence $\bar\eta\circ g'\mid_C=\phi$.

Now let $\dim\nu X\le n$. Let
$q':A\to\R^{n+1}$ be a proper asymptotically Lipschitz map of a
closed subset $A$ of $X$. By Lemma 4.4 there is an asymptotically Lipschitz extension
$q:W\to\R^{n+1}$ over a closed neighborhood $W$ of $A$ with the
property $\|q(a)\|\le\lambda d(a,X\setminus W)+s$ for some numbers $\lambda, s$.
Since $q$ is asymptotically Lipschitz there is an extension to the Higson compactifications
$\bar q:\bar W\to\overline{R^{n+1}}$. Let $\xi:\overline{\R^{n+1}}\to 
\R^{n+1}\cup S^n=B^{n+1}$ be the extension over the Higson compactification
of a radial homeomorphism $h:\R^{n+1}\to Int(B^{n+1})$
given by the formula: $h(t,\theta)=(\frac{t}{t+1},\theta)$.
We note that $\xi(\nu\R^{n+1})\subset\partial B^{n+1}$. 
Since $\dim\nu X\le n$, there is an extension
$\psi:\nu X\to S^n$ of the map $\xi\circ\bar q\mid_{\nu W}:\nu W\to S^n$.
Let $\bar\psi:X\to B^{n+1}$ be an extension of $\psi\cup h\circ q$.
 Consider $g:W\to\R_+$, defined as $g(x)=1+\|q(x)\|$.
Note that $g$ is asymptotically Lipschitz. Denote the constants by 
$\bar\lambda$ and $\bar s$.

 Let $c_n=n\bar\lambda+\bar s$.
Consider $f_n(x)=\frac{c_n}{Var_n\bar\psi(x)}+c_n+1$ for $x\in X\setminus W$
and $f_n9x)=\frac{c_n}{Var_n(\bar\psi\mid_W)(x)}+c_n+1$ for $x\in W$.

Show that $g\le f_n$ on $W$. Let 

\

$Var_n(\bar\psi\mid_W)(x)=
\|h(q(x))-h(q(y))\|=\|\frac{\|q(x)\|}{\|q(x)\|+1}\frac{q(x)}{\|q(x)\|}-
\frac{\|q(y)\|}{\|q(y)\|+1}\frac{q(y)}{\|q(y)\|}\|$. 

\

Denote $a=\|q(x)\|$.
Since $d(x,y)\le n$, we have $a-c_n\le\|q(y)\|\le a+c_n$. Then

\

$Var_n(\bar\psi\mid_W)(x)\le\|(\frac{\|q(x)\|}{\|q(x)\|+1}-\frac{\|q(y)\|}{\|q(y)\|+1})
\frac{q(x)}{\|q(x)\|}\|+\|\frac{\|q(y)\|}{\|q(y)\|+1}(\frac{q(x)}{\|q(x)\|}-
\frac{q(y)}{\|q(y)\|})\|$

\

$\le \frac{c_n}{(a+1)(a-c_n+1)}+(1-\frac{1}{a-c_n+1})
\frac{c_n}{a-c_n}=\frac{c_n}{a-c_n}-\frac{c_n(c_n+1)}{(a+1)(a-c_n+1)(a-c_n)}
\le\frac{c_n}{a-c_n}$

\

provided $a\ge c_n$. 
Therefore $\frac{c_n}{Var_n(\bar\psi\mid_W)(x)}\ge a-c_n$ if $a-c_n\ge 0$.
Hence the inequality $\frac{c_n}{Var_n(\bar\psi\mid_W)(x)}\ge a-c_n$ always holds.
Hence, $f_n(x)=\frac{c_n}{Var_n\bar\psi(x)}+c_n+1\ge a+1=g(x)$.
Let $f$ be a coarsely proper function defined by Proposition 7.5 for $\{f_n\}$ and 
$g$.

Apply Lemma 4.6 to obtain an asymptotically Lipschitz function $\bar g:X\to\R_+$
with $\bar g\le f$ and with $\bar g\mid_{A}=g\mid_{A}$. 
Let $\tilde\lambda$ and $\tilde s$ be its Lipschitz constants.
We define a map
$\tilde q:X\to\R^{n+1}$ as $\tilde q(x)=\bar\psi(x)\bar g(x)$. Note that 
if $x\in A$,
then $\tilde q(x)=\frac{\|q'(x)\|}{\|q'(x)\|+1}\frac{q'(x)}{\|q'(x)\|}(\|q'(x)\|+1)=
q'(x)$. So, $\tilde q$ is an extension of $q'$. Show that $\tilde q$ is asymptotically
Lipschitz. Let $x,y\in X$ be given points. Let $n-1\le d(x,y)\le n$.

Then 

\

$\|\tilde q(y)-\tilde q(x)\|\le
\bar g(y)\|\bar\psi(y)-\bar\psi(x)+\bar\psi(x)-\frac{\bar g(x)}{\bar g(y)}
\bar\psi(y)\|
\le \bar g(y)\|\bar\psi(y)-\bar\psi(x)\|+\frac{|\bar g(y)-\bar g(x)|}
{\bar g(y)}\|\bar\psi(y)\|$

\

$\le \bar g(y)\|\bar\psi(y)-\bar\psi(x)\|+|\bar g(y)-\bar g(x)|$,
since $1\le \bar g\le f$ and $\|\bar\psi(y)\|\le 1$.

\

If $y\in A_n$, then $\bar g(y)\le f(y)\le n\le d(x,y)+1$ and we obtain

\

$\|\tilde q(y)-\tilde q(x)\|\le 2d(x,y)+2+\tilde\lambda d(x,y)+\tilde s$.

\

If $y\in X\setminus A_n$, then $\bar g(y)\le f_n(y)=
\frac{n\bar\lambda+\bar s}{Var_n\bar\psi(y)}+n\bar\lambda+\bar s$ for 
$y\notin W$. Then

\

$\|\tilde q(y)-\tilde q(x)\|\le
\frac{\|\bar\psi(y)-\bar\psi(x)\|}{Var_n\bar\psi(y)}(n\bar\lambda+\bar s)+
2(n\bar\lambda+\bar s)+\tilde\lambda d(x,z)+\tilde s\le$

\

since $\frac{\|\bar\psi(y)-\bar\psi(x)\|}{Var_n\bar\psi(y)}\le 1$, we can
continue

\

$\le 3(n\bar\lambda+\bar s)+\tilde\lambda d(x,y)+\tilde s\le
3((d(x,y)+1)\bar\lambda+\bar s)+\tilde\lambda d(x,y)+\tilde s=
(3\bar\lambda+\tilde\lambda)d(x,y)+(3\bar\lambda+3\bar s+\tilde s)$.

If $y\in W$ we may assume that $x\in W$ as well. Otherwise we apply the
above argument to $x$ instead of $y$. Then
$\frac{\|\bar\psi(y)-\bar\psi(x)\|}{Var_n(\bar\psi\mid_W)(y)}\le 1$ and
the above inequality holds.

\qed

\enddemo

\head \S8 On Coarse Novikov conjectures \endhead

The following coarse statements imply one or another version of
the Novikov conjecture.

{\bf 1. Gromov's conjecture.} {\it A uniformly contractible manifold of
bounded geometry (UC \& BG) is (rationally) hypereuclidean (hyperspherical).}

{\bf 2. Weinberger's conjecture.} {\it The homomorphism $\delta:H^i(\nu X;\Q)\to H^{i+1}_c(X;\Q)$
form the exact sequence of the pair $(\bar X,\nu X)$ is an epimorphism.}

{\bf 3. Coarse Baum-Connes conjecture.} {\it Let $X$ be a proper metric space with
$X\in UC\& BG$. Then the Roe index map
$\mu :K_*^{lf}(X)\to K_*(C(X))$ is an isomorphism (rational isomorphism or rational
monomorphism).}

These conjectures are closely related. The connection between 1 and 2 is discussed 
in [D-F] and the connection between 2 and 3 in [Ro1],[Ro2].

An open $n$-manifold $X$ is called (rationally) {\it hypereuclidean} if there is a morphism
$f:X\to\R^n$ of degree one (nonzero). It is called (rationally) {\it hyperspherical}
if for every $R>0$ there is a short `proper´ map $f:X\to S^n(R)$ onto the standard
sphere of radius $R$ with the degree one (nonzero) [G2]. Here 
`proper´ means that the
complement to some compact set in $X$ goes to a single point.

The most general result supporting the above conjectures is the following
\proclaim{Yu's Theorem [Y2]} If $X\in UC\& BG$ and $X$ admits a coarsely uniform embedding
in the Hilbert space, then the Coarse Baum-Connes conjecture holds for $X$.
\endproclaim

We note that a coarsely uniform embedding is an embedding in the asymptotic category
with Roe's morphisms. Thus, $f:X\to Y$ is a coarsely uniform embedding if there are
two functions $\rho_1,\rho_2:\R_+\to\R_+$ tending to infinity such that
$\rho_1(d_X(x,y))\le d_Y(f(x),f(y))\le\rho_2(d_X(x,y))$. For a geodesic metric space $X$
the image $f(X)$ is not necessarily a totally geodesic subspace, then it means that the
inverse map is not necessarily coarsely Lipschitz.

This theorem in particular implies the previous Yu's result saying that $X$ with
$as\dim X<\infty$ and $X\in UC\& BG$ satisfies the coarse Baum-Connes conjecture.
Here we present some results in that direction.

A sequence of points $\{x_n\}$ with a metric $d$ such that
$\lim_{n\to\infty}d(x_n,\{x_1,\dots, x_{n-1}\})=\infty$ we call a
{\it 0-dimensional asymptotic polyhedron}. We note that such space is 
0-dimensional in any (asymptotic) sense. A 0-dimensional skeleton of an
asymptotic simplex is a 0-dimensional asymptotic polyhedron.
The following lemma first was proved by J. Roe (unpublished).
\proclaim{Lemma 8.1}
The coarse Baum-Connes conjecture holds for 0-dimensional asymptotic
polyhedra.
\endproclaim
\proclaim{Lemma 8.2}
The coarse Baum-Connes conjecture holds for all finite dimensional
asymptotic polyhedra.
\endproclaim
\demo{Proof}
By induction on dimension of asymptotic polyhedra.
The induction starts by virtue of Lemma 8.1.
Let $K$ be an $n$-dimensional asymptotic polyhedron.
For every simplex $\Delta$ we denote by $\frac{1}{2}\Delta\subset\Delta$ an 
$\frac{1}{2}$-homothetic image of $\Delta$ with the common center.
Let $B$ be the union of $\frac{1}{2}\Delta$ for all $n$-simplices in $K$.
Let $A=\overline{K\setminus B}$. Note that $A$ is homotopy equivalent in
a coarse sense to the $n-1$-dimensional skeleton $K^{n-1}$. The space $B$
is homotopy equivalent to an asymptotic 0-dimensional polyhedron. Note that
$C=A\cap B$ is an $n-1$-dimensional asymptotic polyhedron. 
Since the validity of the coarse Baum-Connes conjecture is a coarse homotopy
invariant [Ro2], by the induction assumption
 for spaces 
$A$, $B$ and $C$ the coarse Baum-Connes conjecture holds.
Then by the Five Lemma and the Mayer-Vietoris sequence it holds for $K=A\cup B$.
\qed
\enddemo
\proclaim{Theorem 8.3}
Let $X$ be $ANE(\sA)$ with $as\dim_*X<\infty$.
Then the coarse Buam-Connes conjecture holds for $X$.
\endproclaim
The proof of this theorem follows from Lemma 8.2 and next two lemmas.
\proclaim{Lemma 8.4}
Let $X$ be $ANE(\sA)$ with $as\dim_*X<\infty$.
Then $X$ is homotopy dominated in $\sA$ by an asymptotic polyhedron
of finite dimension.
\endproclaim
\demo{Proof}
Since $X$ is $ANE$ the is a proper map $\alpha:X\to\R_+$ such that
every two $\alpha$-close morphisms to $X$ are homotopic. Let $W$ be a
neighborhood of $X$ in $P(X)$ that admits a retraction $r:W\to X$ with a Lipscitz constant $\lambda$.
Let $d(x)=\frac{1}{\lambda}d_{P(X)}(x,W\setminus r^{-1}(B_{\alpha(x)}(x)))$
 We take
an approximation $\phi:X\to K_d$ of $X$ by an asymptotic polyhedron 
$K_d\subset W$ of
a finite dimension such that $\|\phi\|_A(x)<d(x)$.
Then $r\circ\phi$ is $\alpha$-close to $1_X$.\qed
\enddemo
\proclaim{Lemma 8.5}
Let $Y$ homotopy dominate $X$ in $\sA$ and assume that the coarse Baum-Connes
holds for $Y$. Then it holds for $X$.
\endproclaim
\demo{Proof}
Consider the diagram:
$$
\CD
KX_*(X) @>A_X>> K_*(C^*(X))\\
@Vi_*VV  @Vi'_*VV\\
KX_*(Y) @>A_Y>> K_*(C^*(Y))\\
@Vr_*VV  @Vr'_*VV\\
KX_*(X) @>A_X>> K_*(C^*(X))\\
\endCD
$$
Here homomorphisms $i_*,i'_*,r_*$ and $r'_*$ are induced by morphisms
$i:X\to Y$ and $r:Y\to X$ such that $r\circ i$ is homotopic to the identity map 
$1_X$. Then $i_*$ is a monomorphism and $r'_*$ is an epimorphism.
Since $i_*$ and $A_Y$ are monomorphisms, $A_X$ is a monomorphism.
Since $r'_*$ and $A_Y$ are epimorphic, $A_X$ is an epimorphism.\qed
\enddemo

Let $q:X\to\R_+$ be a function, we say that two maps $\psi_0,\psi_1:Z\to X$
are $q$-close if $d(\psi_0(z),\psi_1(z))<q(\psi_1(z))$. We say that
$\psi_0,\psi_1$ are $q$-homotopic if there is a homotopy $H:Z\times I\to X$ joining 
them such that $diam(H(\{z\}\times I))<q(\psi_0(z))$.
\proclaim{Proposition 8.6}
Let $q:X\to\R_+$ be a function such that $q(x)\le \frac{1}{2}\|x\|$. Then any two 
$q$-close proper maps are proper homotopic.
\endproclaim
\demo{Proof}
Let $H:Z\times I\to X$ be a $q$-homotopy. Let $B_{\rho}(x_0)$ be the ball
of radius $\rho$. Show that $H^{-1}(B_{\rho}(x_0))$ is compact.
Since $\psi_0=H\mid_{Z\times\{0\}}$ is proper, $\psi_0^{-1}(B_{2\rho}(x_0))$ is
compact. We show that  $H^{-1}(B_{\rho}(x_0))\subset
\psi_0^{-1}(B_{2\rho}(x_0))\times I$. For any 
$z\in Z\setminus\psi_0^{-1}(B_{2\rho}(x_0))$ we have
$diam(H(\{z\}\times I))<q(\psi_0(z))<\frac{1}{2}\|\psi_0(z)\|=
\frac{1}{2}d(\psi_0(z),x_0)$.
Hence $d(x_0,H(\{z\}\times I)\ge d(x_0,\psi_0(z))-diam(H(\{z\}\times I)\ge
d(x_0,\psi_0(z))-\frac{1}{2}d(\psi_0(z),x_0)=\frac{1}{2}(\psi_0(z),x_0)>\rho$.
\qed
\enddemo
\proclaim{Theorem 8.7}
Let $X$ be uniformly contractible proper metric space of bounded geometry with $\dim X<\infty$ and 
$as\dim_*X<\infty$.
Then the monomorphism version of the coarse Baum-Connes conjecture
holds true.
\endproclaim
\demo{Proof}
We construct  a morphism $\phi:X\to K$ to a finite dimensional
asymptotic polyhedron $K$ and a map $r:K\to X$ such that 
the composition $r\circ \phi$ is proper homotopic to the identity $1_X$.
Then $r_*\circ\phi_*=id$. 
Since for UC and BG spaces $KX_*=K_*^{lf}$, the diagram of Lemma 8.5
turns into
$$
\CD
K_*(X) @>A_X>> K_*(C^*(X))\\
@V\phi_*VV  @V\phi'_*VV\\
K_*(Y) @>A_Y>> K_*(C^*(Y))\\
@Vr_*VV  @.\\
K_*(X) @.\\
\endCD
$$
Then Lemma 8.2 and the argument of Lemma 8.5 would imply that $A_X$ is
a monomorphism. 

Let $\dim X=m$. Since $X\in UC$ there is a proper monotone function
$g_m:\R_+\to\R_+$ such that for any function $q:X\to\R_+$ any two $q$-close
maps of $m$-dimensional space to $X$ are $g_m\circ q$-homotopic.

Let $as\dim_*X=n$. Let $S:\R_+\to\R_+$ be a contractibility function of $X$.
We define $T(t)=\frac{1}{2}S^{-1}(\frac{t}{2})$, $\bar f=T^n=T\circ\dots\circ T$ and
define $f:X\to\R_+$ as $f(x)=\bar f(\frac{1}{2}g_m^{-1}(\frac{1}{4}\|x\|))$. 
Let $\phi:X\to K_f$ be a short
map to an asymptotic polyhedron with $\|\phi\|_U<\frac{1}{4}f$. Then we
can subdivide $K_f$ to obtain an asymptotic polyhedron with the property that
$diam(\phi^{-1}(\Delta))<f(x)$ for all $x\in\phi^{-1}(\Delta)$. We can achieve
this for $\Delta$ lying outside of some compact set. We may (and we will) ignore
this restriction in the further argument.

By induction we construct a map $r:K_f^{(k)}\to X$ of $k$-dimensional skeleton.
For any $v\in K_f^{(0)}$ we define $r(v)\in\phi^{-1}(v)$. Since for every
1-simplex $[u,v]$ with vertices $u,v$, $diam\phi^{-1}([u,v])<f(r(u))\le
\bar f(\frac{1}{2}g_m^{-1}(\frac{1}{4}\|r(u)\|))=T^n(\frac{1}{2}g_m^{-1}(\frac{1}{4}\|r(u)\|))=\frac{1}{2}S^{-1}
(\frac{1}{2}T^{n-1}(\frac{1}{2}g_m^{-1}(\frac{1}{4}\|r(u)\|)))$, 
there is an extension of $r$
over $[u,v]$ with $diam(r([u,v]))<\frac{1}{2}T^{n-1}(\frac{1}{2}g_m^{-1}(\frac{1}{4}\|r(u)\|))$. 
We extend $r$ over $K^{(1)}_f$ in this manner.
Note that for any 2-simplex
$\sigma$ spanned by $u,v,w$, $diam(r(\sigma^{(1)}))<T^{n-1}(\frac{1}{2}g_m^{-1}(\frac{1}{4}\|r(u)\|))$. Therefore
there is an extension of $r$ over $\sigma$ with $diam(r(\sigma))<\frac{1}{2}
T^{n-2}(\frac{1}{4}\|r(u)\|)$.
Finally we will construct a map $r:K^{(n)}_f=K\to X$ with
$diam(r(\Delta))<\frac{1}{2}g_m^{-1}(\frac{1}{4}\|r(u)\|)$ for any $n$-simplex 
$\Delta$ where $u$ is a vertex 
in $\Delta$. Then for any point $x\in X$ the distance
$d(x,r\phi(x))$ can be estimated from above as $d(x,r(u))+d(r(u),r\phi(x))<
f(x)+\frac{1}{2}g_m^{-1}(\frac{1}{4}\|r(u)\|)<f(x)+\frac{1}{2}g_m^{-1}(
\frac{1}{4}(\|x\|+f(x)))
\le \frac{1}{2}g_m^{-1}(\frac{1}{4}\|x\|)+\frac{1}{2}g_m^{-1}(\frac{1}{4}\|x\|
+\frac{1}{8}\|x\|)
\le g_m^{-1}(\frac{1}{2}\|x\|)$. Here we used the fact
that $T(t)<t$ and $g_m^{-1}(t)\le t$. Hence $1_X$ and $r\circ\phi$ are
$\frac{1}{2}\|x\|$-homotopic. Proposition 8.6 implies that $1_X$ and $r\circ\phi$
are proper homotopic.\qed
\enddemo
\proclaim{Theorem 8.8}
Let $X$ be a uniformly contractible $n$-manifold with asymptotic dimension
$as\dim X=n$. Then $X$ is hyperspherical.
\endproclaim
\demo{Proof}
Let $\Delta$ be the standard simplex of dimension $n$. Let $c_n$ be the
Lipschitz constant of a map $\nu$ taking $\Delta$ to a unit hemisphere $S_+$ 
homeomorphically except 
 the complement in $\partial\Delta$ to one $n-1$-dimensional face goes to 
 a point. Let $\epsilon$ be given. Consider a short map $\phi:X\to K$ to
 $n$-dimensional polyhedron with mesh $m> \frac{c_n}{\epsilon}$. We consider a
 sphere $S$ of sufficiently large radius $R$ centered at $s_0=\phi(x_0)$. Clearly,
 $S$ separates $s_0$ and $\infty$. Take a smallest subcomplex $N\subset K$
 that contains $S$. For big enough $R$ the complex $N$ separates $s_0$ and $\infty$.
 Since $\dim K=n$, $n-1$-dimensional skeleton $N^{(n-1)}$ separates $s_0$ and
 $\infty$. Let $N'\subset N^{(n-1)}$ be the boundary of the component containing
 $s_0$.
 Then $M=\phi^{-1}(N')$ separates $x_0$ and $\infty$. We may assume that
 $\phi$ is a light simplicial map with respect to some triangulation of $X$
 and a small subdivision of $K$. Then $M$ is a polyhedron.
 The boundary of the component containing $x_0$ defines a cycle $c$ in $M$
 that generates the $n-1$-dimensional homology group of $X\setminus\{x_0\}$.
 Let $C=supp(c)$. Since $X\in UC$, there is an 'approximate' lift
 $\alpha:K\to X$ of $\phi$. If $R$ is big enough then 
 $\alpha\circ\phi\mid_C$ is homotopic to $1_C$ in $X\setminus\{x_0\}$.
 Then $(\phi\mid_C)_*(c)\ne 0$. Hence there is an $(n-1)$-simplex $\sigma$ in
 $N'$ such that $\phi\mid_{\phi^{-1}(\sigma)}$ is essential.
 An easy diagram chasing
 shows that the degree of the homomorphism
 $\Bbb Z=H_{n-1}(C)\to H_{n-1}(\phi^{-1}(N'))\to H_{n-1}(\phi^{-1}(N'),
 \phi^{-1}(N'\setminus Int(\sigma))\to H_{n-1}(N',N'\setminus Int(\sigma))=
H_{n-1}(\sigma,\partial\sigma)=\Bbb Z$ is one.
 Let $W$ be the star of $\sigma$ in $K$, i.e. $W$ is the union of all
 simplices containing $\sigma$. We define a $\frac{c_n}{m}$-Lipschitz map
 $\psi:W\to S^{n}$ to the unit sphere using the map $\nu$ in
 such a way that all $n$-simplices $\Delta\supset\sigma$ lying in the component
 of $s_0$ go to the lower hemisphere and all other simplices 
 $\Delta\supset\sigma$ are mapped
 to upper hemisphere. Then the map $\bar\psi:X\to S^n$ defined as
 $\psi\circ\phi$ on $\phi^{-1}(W)$ and as a constant map on the rest of $X$
 is $\epsilon$-Lipschitz of degree one.
\enddemo
\proclaim{Theorem 8.9}
Let $X$ be a metric space with slow dimension growth,
then the Coarse Baum-Connes conjecture holds for $X$.
\endproclaim
The proof is based on Yu's theorem and the following
Lemma which is due to Higson-Roe and Yu [H-R2],[Y2].
A probability measure $\mu$ on a discrete space $Z$
can be treated as a function $\mu:Z\to\R_+$ with
$\Sigma_{z\in Z}\mu(z)=1$, the support $supp(\mu)$
consists of those $z\in Z$ where $\mu$ is nonzero.
By $\|\ \|_1$ we denote the $l_1$-norm on a space of functions.

\proclaim{Lemma 8.10 [H-R2],[Y2]}
Let $Z$ be a discrete space of bounded geometry and assume
that there is a sequence of maps $a^n:Z\to P(Z)$
to probability measures such that

(1) For any $n$ there is $R$ such that 
$supp(a^n(z))\subset\{z'\in Z\mid d(z,z')<R\}$.

(2) For any $K>0$, $\lim_{n\to\infty}\sup_{d(z,w)<K}
\|a^n(z)-a^n(w)\|_1=0$.

Then $Z$ admits a coarsely uniform embedding into
a Hilbert space.
\endproclaim
\demo{Proof of Theorem 8.9}
It suffices to embed an $\epsilon$-dense discrete
subset $Z\subset X$. For given $L$ there is an
open cover $\sU$ of $Z$ with the Lebesgue number
$L$ of multiplicity $d(L)$. Let $\sU=\{U_i\}_{i\in J}$.
We define $\phi_i(x)=\frac{d(x,Z\setminus U_i)}{
\Sigma_{i\in J}d(x,Z\setminus U_j)}$.
Denote by $\Sigma_x=\Sigma_{i\in J}d(x,Z\setminus U_j)$.
For every $i$ we fix $y_i\in U_i$.
Define $a^L(z)=\Sigma_{i\in J}\phi_i(z)\delta_{y_i}\in
P(Z)$. Note that $supp(a^L(z))=\{y_i\mid\phi_i(z)\ne 0\}=
\{y_i\mid z\in U_i\}\subset\{z'\in Z\mid d(z,z')<R\}$ where 
$mesh\sU<R$. 

To check (2) note that

$\|a^L(z)-a^L(w)\|_1=\|\Sigma_{i\in J}\phi_i(z)\delta_{y_i}-
\Sigma_{i\in J}\phi_i(w)\delta_{y_i}\|_1=
\|\Sigma_{i\in J}(\phi_i(z)-\phi_i(w))
\delta_{y_i}\|_1$

\

$\le\Sigma_J|\phi_i(z)-\phi_i(w)|=
\Sigma_J|\frac{d(z,Z\setminus U_i)}{\Sigma_z}-
\frac{d(w,Z\setminus U_i)}{\Sigma_w}|\le
\Sigma_J|\frac{d(z,Z\setminus U_i)}{\Sigma_z}-
\frac{d(w,Z\setminus U_i)}{\Sigma_z}|+$

\

$\Sigma_Jd(w,Z\setminus U_i)|\frac{1}{\Sigma_z}-\frac{1}{\Sigma_w}|
\le\frac{1}{\Sigma_z}\Sigma_J|d(z,Z\setminus U_i)-d(w,Z\setminus U_i)|
+\Sigma_w|\frac{1}{\Sigma_z}-\frac{1}{\Sigma_w}|
\le$ 

\ 

$\frac{2d(L)}{\Sigma_z}d(z,w)+
\frac{|\Sigma_w-\Sigma_z|}{\Sigma_z}\le
\frac{2d(L)}{L}d(z,w)+\frac{2d(L)}{L}d(z,w)\le
4K\frac{d(L)}{L} \to 0$ as $L\to\infty$.
\qed
\enddemo

In classical dimension theory of infinite dimensional spaces there is
a special class of spaces having the Property C. Compact metric spaces
with the Property C are  having some features of finite dimensional spaces
(see [vM-M]). Below we define an asymptotical analog of the Property C.
\proclaim{Definition}
A metric space $X\in BG$  has an asymptotic property C if for
every sequence of numbers $R_1\le R_2\le R_3\dots$ there exists a finite
sequence of uniformly bounded families $\{\sU_i\}^k_{i=1}$ of 
open subsets of $X$ such that the union $\cup_{i=1}^k\sU_i$ is a cover of $X$
and every family $\sU_i$ is $R_i$-disjoint.
\endproclaim
The following theorem is close to Theorem 8.9.
\proclaim{Theorem 8.11} The Coarse Baum-Connes conjecture holds for
metric spaces with the asymptotic property C.
\endproclaim
\demo{Proof}
We apply Lemma 8.10 and the Theorem of Yu. For every $n$ we define
a map $a^n:X\to P(X)$ as follows. First we assume that $X$ is a discrete space.
Let $R_i=n^i$ and let $\sU_1,\dots,\sU_k$ be a sequence of $2R_i$-disjoint
uniformly bounded families from the definition of the asymptotic property C.
For every $U^i_j\in\sU_i$ we fix a point $x^i_j\in U^i_j$. We define

$$
\phi^i_j(x)=\cases\frac{R_i}{2},&if \ x\in U^i_j;\\
\max\{0,\frac{R_i}{2}-d(x,U^i_j)\},&otherwise
\endcases
$$

and define $b^n(x)=\Sigma_{i,j}n^{k-i+1}\phi^i_j(x)\delta_{x^i_j}$.
We set $a^n(x)=\frac{b^n(x)}{\|b^n(x)\|_1}$.
Note that for every $x\in X$ and any $i$ there is at most one $j=j_x(i)$
such that $\phi^i_j(x)\ne 0$. Then $supp(a^n(x)$ is contained in the ball $B_r(x)$ with $r=\max\{diam\sU_i+2R_i\}$. Thus the condition (1) of
Lemma 8.10 holds.

In order to verify the condition (2) first we show that
$\|b^n_x\|_1\ge\frac{n^{k+1}}{2}$ for all $x$. Indeed,
$\|b^n_x\|_1=\Sigma_in^{k-i+1}|\phi^i_{j(i)}(x)|\ge n^{k-i+1}\frac{R_i}{2}=
\frac{n^k+1}{2}$ where $i$ is chosen with $x\in U^i_{j(i)}$.

Without loss of generality we may assume that $\|b^n_y\|_1\ge\|b^n_x\|_1$.
Then $\|b^n_x-\|b^n_x\|_1a^n_y\|_1\le\|\frac{\|b^n_x\|_1}{\|b^n_y\|_1}b^n_y-b^n_y\|_1+\|b^n_y-b^n_x\|_1$

\

$=\|b^n_y\|_1\frac{|\|b^n_x\|_1-\|b^n_y\|_1|}{\|b^n_y\|_1}+\|b^n_y-b^n_x\|_1\le
2\|b^n_y-b^n_x\|_1$

\

$2\|\Sigma_in^{k-i+1}\phi^i_{j_x(i)}(x)\delta_{x^i_{j_x(i)}}-
\Sigma_in^{k-i+1}\phi^i_{j_y(i)}(y)\delta_{x^i_{j_y(i)}}\|_1\le$

\

$2\|\Sigma_{i,j(i)=j_x(i)=j_y(i)}n^{k-i+1}(\phi^i_{j(i)}(x)-\phi^i_{j(i)}(y))
\delta_{x^i_{j(i)}}\|_1+$

(here we consider $j(i)$ equal to any of the indexes $j_x(i)$ or $j_y(i)$
when one of them is not defined)

\

$2\|\Sigma_{i,j_x(i)\ne j_y(i)}n^{k-i+1}\phi^i_{j_x(i)}(x)\delta_{x^i_{j_x(i)}}\|_1+
2\|\Sigma_{i,j_x(i)\ne j_y(i)}n^{k-i+1}\phi^i_{j_y(i)}(y)\delta_{x^i_{j_y(i)}}\|_1$

\

$\le 2\Sigma_in^{k-i+1}|\phi^i_{j(i)}-\phi^i_{j(i)}(y)|+4\Sigma_in^{k-i+1}d(x,y)$.

Here we used that for $j_x(i)\ne j_y(i)$ we have $d(x,y)\ge R_i\ge 2\phi^i_j$.

Note that $|\phi^i_{j(i)}-\phi^i_{j(i)}(y)|\le d(x,y)$. Then

$\Sigma_in^{k-i+1}|\phi^i_{j(i)}-\phi^i_{j(i)}(y)|\le
\Sigma_{i=1}^kn^{k-i+1}d(x,y)=n\frac{n^k-1}{n-1}d(x,y)$.

Then $\|b^n_x-\|b^n_x\|_1a^n_y\|_1\le 6n\frac{n^k-1}{n-1}d(x,y)$.
Therefore $\|a^n_x-a^n_y\|_1=\frac{1}{\|b^n_x\|}\|b^n_x-\|b^n_x\|_1a^n_y\|_1\le$

\

$\frac{6n}{\|b^n_x\|}\frac{n^k-1}{n-1}d(x,y)\le\frac{12n}{n^{k+1}}\frac{n^k-1}{n-1}d(x,y)\le \frac{12}{n-1}d(x,y)$. 

If $d(x,y)\le K$ we have that
$\lim_{n\to\infty}frac{12}{n-1}d(x,y)=0$.\qed
\enddemo

\head \S9 Principle of Descent and the Higson corona \endhead

The Descent Principle is the statement that the original Novikov higher
signature conjecture for geometrically finite group can be derived from
its coarse counterpart. The main coarse analog of the Novikov Conjecture
is the coarse Baum-Connes conjecture which was considered in the previous
section. Here we consider coarser statements which also imply the
Novikov Conjecture.  Let $X$ be a universal cover
of finite aspherical polyhedron $B\Gamma$ supplied with a metric lifted from
$B\Gamma$. Then each of the following four conditions implies the 
Novikov conjecture for $\Gamma$.

\proclaim{(CPI) \cite{C-P}} There is an equivariant rationally acyclic metrizable
compactification $\hat X$ of $X$ such that the action of $\Gamma$ is small
at infinity.
\endproclaim
\proclaim{(CPII) \cite{C-P2}} There is an equivariant rationally acyclic 
(possibly nonmetrizable) compactification $\hat X$ of $X$ with a system of 
covers ${\alpha}$ of $Y=\hat X\setminus X$ by boundedly saturated sets
such that the projection to the inverse limits of the nerves of $\alpha$ 
induces an isomorphism $H_*(Y;\Q)\to H_*(\lim_{\leftarrow}N(\alpha);\Q)$.
 \endproclaim
\proclaim{(FW) \cite{F-W},\cite{D-F}} The boundary homomorphism
$\delta: H^{n-1}(\nu X;\Q)\to H^n_c(X;\Q)$ is an equivariant split surjection.
\endproclaim
\proclaim{(HR) \cite{Ro1}} There is an equivariant rationally acyclic Higson
dominated compactification $\hat X$ of $X$.
\endproclaim
An action of $\Gamma$ is
{\it small at infinity} for a given compactification $\bar X$ of $X$ 
if for every $x\in\bar X\setminus X$ and a neighborhood $U$ of $x$ in $\bar X$,
for every compact set $C\subset X$ there is a smaller neighborhood $V$ such that
$g(C)\cap V\ne\emptyset$ implies $g(C)\subset U$ for all $g\in\Gamma$.

An open set $U\subset Y=\hat X\setminus X$ is called 
{\it boundedly saturated} if for for every closed set $C\subset\hat X$ with
$C\cap Y\subset U$ the closure of any $r$-neighborhood $N_r(C\cap X)$
satisfies $\overline{N_r(C\cap X)}\cap Y\subset U$.

We consider a $\Gamma$-invariant metric on $X$. Since $B\Gamma$ is a finite
complex, the Higson corona of $X$ does not depend on choice of metric and coincides
with the Higson corona of $\Gamma$.

\proclaim{Proposition 9.1 \cite{D-F}}
The action of $\Gamma$ on $X$ is small at infinity
for a compactification $\bar X$ if and only if $\bar X$ is Higson dominated.
\endproclaim

Let $\Gamma$ be a finitely generated group with the word metric 
$d(a,b)=l(a^{-1}b)$ where $l(w)$ denotes the minimal length of a word presenting
 $w\in\Gamma$ with a fixed finite set of generators.
A group $\Gamma$ acts on itself by isometries: $g:\Gamma\to\Gamma$ by the
formula $g(a)=ga$. This action is called the {\it left action}.
\proclaim{Proposition 9.2}
The left action of $\Gamma$ on itself can be extended to an action on the
Higson compactification $\bar\Gamma$.
\endproclaim
\demo{Proof}
An isometry $g:\Gamma\to\Gamma$ is extendible over the Higson corona to
a continuous map $g:\bar\Gamma\to\bar\Gamma$. This implies that the whole
action of $\Gamma$ is extendible.
\enddemo
\proclaim{Lemma 9.3}
Let $X$ be as above.
The following are equivalent.
\roster
\item{}The boundary homomorphism
$\delta: H^{n-1}(\nu X;\Q)\to H^n_c(X;\Q)$ is an equivariant split surjection.
\item{} There is a Higson dominated equivariant compactification $\hat X$ of
$X$ such that the boundary homomorphism
$\hat\delta: H^{n-1}(\hat X\setminus X;\Q)\to H^n_c(X;\Q)$ is an 
equivariant split surjection.
\item{}There is a metrizable Higson dominated equivariant compactification $\hat X$ of
$X$ such that the boundary homomorphism
$\hat\delta: H^{n-1}(\hat X\setminus X;\Q)\to H^n_c(X;\Q)$ is an 
equivariant split surjection.
\item{}There is an equivariant metrizable Higson dominated
compactification $\hat X$ of $X$ such that the boundary homomorphism
$H^{lf}_*(X;\Q)\to H^s_{*-1}(\hat  X\setminus X;\Q)$ is an equivariant 
split injection.
\endroster
\endproclaim
Here $H^s_*$ stands for the Steenrod homology.
\demo{Proof}
1)$\Rightarrow$ 2). Because of Proposition 9.2 this implication follows.

2)$\Rightarrow$ 1). Let $\xi^*:H^{n-1}(\hat X\setminus X;\Q)\to H^{n-1}(\nu X;\Q)$
be a homomorphism generated by the domination $\xi:\bar X\to\hat X$. The map
$\xi$ is equivariant, since it is equivariant on a dense subset $X$. Let
$s':H^n_c(X;\Q)\to H^{n-1}(\hat X\setminus X;\Q)$ be an equivariant splitting
of the boundary homomorphism. Then $\xi^*\circ s'$ is an equivariant splitting
of $\delta: H^{n-1}(\nu X;\Q)\to H^n_c(X;\Q)$.

3)$\Rightarrow$ 2). Obvious.

3)$\Leftrightarrow$ 4). This is the standard duality (of vector spaces)
between homology and cohomology with coefficients in a field.

For a proof of the implication 2)$\Rightarrow$ 3) we need to develop some
technique.
\enddemo

A {\it directed} set $A$ is a partially ordered set with the property that
for any two elements $\alpha,\beta\in A$ there exists $\gamma\in A$ with
$\gamma\ge\alpha$ and $\gamma\ge\beta$. A directed set is called 
{\it $\sigma$-complete} if for any countable chain $C\subset A$ there is
a supremum $\sup C\in A$. A subset $B\subset A$ of a $\sigma$-complete set $A$
is called {\it $\sigma$-closed} if for any countable chain $C$ in $B$ we have
$\sup C\in B$. A subset $B\subset A$ is called {\it cofinal} if for any 
$\alpha\in A$ there is $\beta\in B$ with $\beta\ge\alpha$.
\proclaim{Proposition 9.4} The intersection of countably many $\sigma$-closed
cofinal subsets $\cap B_i\subset A$ is $\sigma$-closed cofinal.
\endproclaim
\demo{Proof} It is clear that $\cap B_i$ is $\sigma$-closed. First we show that
the intersection of two $\sigma$-closed cofinal sets $B$ and $B'$ is cofinal. For that
for any $\alpha\in A$ we construct a sequence $\beta_i$ such that
1) $\beta_0=\alpha$, $\beta_{2k-1}\in B$ and $\beta_{2k}\in B'$, $k>0$;
2) $\beta_i<\beta_{i+1}$. Then $\sup\{\beta_i\}\in B\cap B'$.

Next, we construct a sequence $\alpha_i$ such that: 1) $\alpha_i\le\alpha_{i+1}$,
2) $\alpha_i\in\cap_k^iB_k$ and $\alpha_0\in A$ is given . Then 
$\sup\{\alpha_i\}\in\cap B_i$.\qed
\enddemo
An inverse (direct) system in a given category $\sC$ over a $\sigma$-complete
ordered set $A$ is called {\it $\sigma$-continuous} if

$\lim_{\leftarrow}\{X_{\alpha_i}\mid \alpha_i\in C\}=X_{\sup C}$
\ \ \ \ ($\lim_{\rightarrow}\{H_{\alpha_i}\mid\alpha_i\in C\}=X_{\sup C}$)

for every countable chain $C$.
\proclaim{Schepin Spectral Theorem}[Sc],[Ch].
Let $\{X_{\alpha}\mid\alpha\in A\}$ and $\{Y_{\alpha}\mid\alpha\in A\}$
be inverse $\sigma$-continuous systems of compact metric spaces.

1) Existence. Let $f:X=\lim\{X_{\alpha}\}\to Y=\lim\{Y_{\alpha}\}$ be a continuous
map. Then there exists a $\sigma$-closed cofinal subset $B\subset A$ and a
morphism of spectra 

$\{f_{\beta}\}_{\beta\in B}:\{X_{\beta}\mid\beta\in B\}\to
\{Y_{\beta}\mid \beta\in B\}$ 
such that $f=\lim_{\leftarrow}\{f_{\beta}\mid\beta\in B\}$.

2) Uniqueness. Let $\{f_{\alpha}\}_{\alpha\in A}$ and 
$\{g_{\alpha}\}_{\alpha\in A}$ be two morphisms between the spectra
$\{X_{\alpha}\mid\alpha\in A\}$ and $\{Y_{\alpha}\mid\alpha\in A\}$ with
$\lim_{\leftarrow}\{f_{\alpha}\}=\lim_{\leftarrow}\{g_{\alpha}\}$. Then
there exists a $\sigma$-closed cofinal subset $B\subset A$ such that
$f_{\beta}=g_{\beta}$ for $\beta\in B$.
\endproclaim
\proclaim{Corollary} For any homeomorphism
$h:X=\lim\{X_{\alpha}\}\to Y=\lim\{Y_{\alpha}\}$ 
there exists a $\sigma$-closed cofinal subset $B\subset A$ and an
isomorphism of spectra $\{h_{\beta}\}_{\beta\in B}:\{X_{\beta}\mid\beta\in B\}\to
\{Y_{\beta}\mid \beta\in B\}$ such that 
$h=\lim_{\leftarrow}\{h_{\beta}\mid\beta\in B\}$.
\endproclaim
\demo{Proof} Apply existence part for $h$ and $h^{-1}$ and then apply
Proposition 9.4.
\enddemo
For every compact Hausdorff space $X$ one can define a continuous
$\sigma$-spectrum $\{X_{\alpha}\mid\alpha\in A\}$ as follows. Consider
an imbedding $X\subset I^B$. Let $A$ be the set of all countable subsets of
$B$, define $X_{\alpha}=\pi_{\alpha}(X)$, where $\pi_{\alpha}:I^B\to I^{\alpha}$
is a projection, $\alpha\in A$. All bonding maps are defined similarly.

There is a dual theorem to the Schepin Spectral theorem for $\sigma$-continuous
direct system of countable CW-complexes [Dr-Dy1]. In this paper we need the dual
theorem in the category of abelian groups.
\proclaim{Dual Spectral Theorem}
Let $\{H_{\alpha}\mid\alpha\in A\}$ and $\{G_{\alpha}\mid\alpha\in A\}$
be direct $\sigma$-continuous systems of countable abelian groups.

1) Existence. Let $f:H=\lim\{H_{\alpha}\}\to G=\lim\{G_{\alpha}\}$ be a 
homomorphism. Then there exists a $\sigma$-closed cofinal subset $B\subset A$ 
and a
morphism of spectra 

$\{f_{\beta}\}_{\beta\in B}:\{H_{\beta}\mid\beta\in B\}\to
\{G_{\beta}\mid \beta\in B\}$ 
such that 
$f=\lim_{\rightarrow}\{f_{\beta}\mid\beta\in B\}$.

2) Uniqueness. Let $\{f_{\alpha}\}_{\alpha\in A}$ and 
$\{g_{\alpha}\}_{\alpha\in A}$ be two morphisms between the spectra
$\{H_{\alpha}\mid\alpha\in A\}$ and $\{G_{\alpha}\mid\alpha\in A\}$ with
$\lim_{\rightarrow}\{f_{\alpha}\}=\lim_{\rightarrow}\{g_{\alpha}\}$. Then
there exists a $\sigma$-closed cofinal subset $B\subset A$ such that
$f_{\beta}=g_{\beta}$ for $\beta\in B$.
\endproclaim
\demo{Proof} Let $p_{\alpha}^{\beta}:H_{\alpha}\to H_{\beta}$ and
$q_{\alpha}^{\beta}:G_{\alpha}\to G_{\beta}$ be the bonding maps, 
$\alpha\le\beta$. Denote 
$p_{\alpha}=\lim_{\beta}p^{\beta}_{\alpha}:H_{\alpha}\to H$ and
$q_{\alpha}=\lim_{\beta}q^{\beta}_{\alpha}:G_{\alpha}\to G$

1). First we show that the set $M_H=\{\alpha\in A\mid Ker(p_{\alpha})=0\}$ is
cofinal in $A$. For any $\alpha_0\in A$ we built the chain
$\alpha_0<\alpha_1<\dots$ such that 
$Ker(p^{\alpha_{i+1}}_{\alpha_i})=Ker(p_{\alpha_i})$. Let 
$\alpha=\sup\{\alpha_i\}$ and let $a\in Ker(p_{\alpha})$. 
 Since the spectrum is $\sigma$-continuous,
$a=p^{\alpha}_{\alpha_i}(c)$ for some $i$ and some $c$. Then $c\in Ker(p_{\alpha_i})$,
therefore $p^{\alpha_{i+1}}_{\alpha_i}(c)=0$. Hence $a=0$.

Let $M=M_H\cap M_G$. Note that $M$ is $\sigma$-closed and cofinal.
We show that the set $B=\{\beta\in M\mid{\text{there exists
a homomorphism}}\ f_{\beta}:H_{\beta}\to G_{\beta}\ 
{\text with}\ f\circ p_{\beta}=q_{\beta}\circ f_{\beta}\}$ is
cofinal in $A$. Let $\alpha\in A$ be given. For any $a\in H_{\alpha}$ there is
$\beta(a)\in A$ such that $fp_{\alpha}(a)=q_{\beta(a)}(b(a))$ for some element
$b(a)\in G_{\beta(a)}$. Enumerate elements of $H_{\alpha}$ as $a_1,a_2,\dots$ and
construct a sequence $\alpha<\beta(a_1)<\beta(a_2)<\dots$. Let $\gamma
=\sup\{\beta(a_i\}$. Then we define a map 
$f_{\alpha\gamma}:H_{\alpha}\to G_{\gamma}$ by the formula
$f_{\alpha\gamma}(a)=q^{\gamma}_{\beta(a)}(b(a))$. Note that $f_{\alpha\gamma}$ is
a homomorphism with $f\circ p_{\alpha}=q_{\gamma}\circ f_{\alpha\gamma}$.
The latter is obvious, we check that $f_{\alpha\gamma}$ is a homomorphism.
Since $q_{\gamma}$ is a monomorphism, the equalities
$q_{\gamma}f_{\alpha\gamma}(a+a')=fp_{\alpha}(a+a')=fp_{\alpha}(a)+
fp_{\alpha}(a')=
q_{\gamma}f_{\alpha\gamma}(a)+q_{\gamma}f_{\alpha\gamma}(a')$ imply that
$f_{\alpha\gamma}(a+a')=f_{\alpha\gamma}(a)+f_{\alpha\gamma}(a')$. Define
$\alpha_1=\gamma$ and repeat the procedure to obtain $\alpha_2$ and so on.
We will obtain a chain $\alpha<\alpha_1<\alpha_2<\dots$ and a sequence of
homomorphisms $f_{\alpha_i\alpha_{i+1}}:H_{\alpha_i}\to G_{\alpha_{i+1}}$
such that all squares are commutative. Then there exists a homomorphism
$f_{\bar\alpha}:H_{\bar\alpha}\to G_{\bar\alpha}$ with $f_{\bar\alpha}=
\lim_{\to}f_{\alpha_i\alpha_{i+1}}$ such that 
$fp_{\bar\alpha}=q_{\bar\alpha}f_{\bar\alpha}$. Clearly $B$ is $\sigma$-closed and
$\{f_{\beta}\}_{\beta\in B}$ is a morphism of spectra.

2). Take $B=M$, then the result follows.

\enddemo

\demo{Proof of 2)$\Rightarrow$ 3) of Lemma 9.3} 
According to the remark after Schepin Spectral Theorem 
we can 
present $\hat X$ as the inverse limit of a $\sigma$-continuous system of compact
metric spaces $\hat X=
\lim_{\leftarrow}\{\hat X_{\alpha};p_{\alpha}^{\beta}\mid\alpha\in A'\}$.
It easy to verify thet the set $A=\{\alpha\mid p_{\alpha}^{-1}p_{\alpha}(X)=X\}$
is $\sigma$-closed cofinal in $A'$.
For any element $\gamma\in\Gamma$ we apply the Schepin Spectral Theorem
to a homeomorphism $\gamma:\hat X\to\hat X$ to obtain a $\sigma$-closed
cofinal subset $B_{\gamma}\subset A$ with corresponding isomorphism of spectra.
By Proposition 9.4 the set $B=\cap_{\gamma\in\Gamma}B_{\gamma}$ is $\sigma$-closed
cofinal. Applying the uniqueness part of the Schepin Spectral Theorem, we
may assume that $\hat X_{\alpha}$ is an equivariant compactification of $X$
such that all bonding maps are equivariant.

The homomorphism $\hat\delta$ is induced by the quotient map
$\hat X\cup cone(\hat X\setminus X)\to\Sigma(\hat X\setminus X)$. So,
$\hat\delta:G=H^n(\Sigma(\hat X\setminus X);\Q)\to H=H^n(\hat X\cup cone(\hat X\setminus X);\Q)$.
The condition 2) says that $\hat\delta$ admits an equivariant splitting 
$s:H\to G$ i.e. $\hat\delta\circ s=id$. By part one of the Dual Spectral Theorem
there exists $B_1\subset B$ such that $s$ is the limit of a morphism
$\{s_{\beta}:H_{\beta}\to G_{\beta}\mid\beta\in B_1\}$. By part two of the Dual
Spectral Theorem we may assume that that $\gamma^{-1}s_{\beta}\gamma=s_{\beta}$
and $\hat\delta_{\beta}s_{\beta}=id$ for all $\beta\in B_1$ and all 
$\gamma\in\Gamma$. Then for any $\beta\in B_1$ a compact $\hat X_{\beta}$
can serve as a metrizable Higson dominated compactification required in 3).

\enddemo

Let $c_1X$ and $c_2X$ be two compactifications of $X$, the maximal compactification
$cX$ dominated by both $c_1X$ and $c_2X$ is called the minimum of $c_1X$ and 
$c_2X$ and is denoted as $\min\{c_1X,c_2X\}$. We recall that compactifications of $X$ are
in one-to-one correspondence with totally bounded uniformities on $X$ [En].
Two sets $A,B\subset X$ have nonempty intersection at infinity, $\bar A\cap\bar B\ne\emptyset$,
if and only if $(A\times B)\cap V\ne\emptyset$ for all entourages $V\in\sU$.
Let $\sU_1$ and $\sU_2$ be two uniformities on $X$ corresponding compactifications
$c_1X$ and $c_2X$, then the uniformity $\sU$ corresponding to $\min\{c_1X,c_2X\}$
is defined as $\sU=\sU_1\cap\sU_2$.

For a compactification $cX$ of locally compact space $X$ we denote by
$\sP(cX)$ the set of pairs $(A,B)$ of closed sets in $X$ with $\bar A\cap\bar 
B\ne\emptyset$. 
\proclaim{Proposition 9.5} Let $X$ be a locally compact space, let $c_1X$ and
$c_2X$ be two compactifications of $X$ and let $cX=\min\{c_1X,c_2X\}$.
Then $\sP(cX)=\sP(c_1X)\cup\sP(c_2X)$.
\endproclaim
\demo{Proof}
It is clear that $\sP(cX)\supset\sP(c_1X)\cup\sP(c_2X)$. Let $(A,B)\in\sP(cX)$.
Assume that $(A,B)\not\in\sP(c_iX)$, $i=1,2$. Then there exist $V_i\in\sU_i$
such that $(A\times B)\cap V_i=\emptyset$. Therefore $(A\times B)\cap(V_1\cup V_2)=\emptyset$.
By the definition of uniformity $V_1\cup V_2\in\sU_i$ for $i=1,2$. Then
$V_1\cup V_2\in\sU$. Contradiction.\qed
\enddemo

The homomorphism
$H_*(Y;\Q)\to H_*(\lim_{\leftarrow}N(\alpha);\Q)$ in $CPII$
is an isomorphism if the system $\{\alpha\}$ is cofinal. We introduce
the condition.
\proclaim{($CPII'$)} There is an equivariant rationally acyclic 
(possibly nonmetrizable) compactification $\hat X$ of $X$ with a 
cofinal system of 
covers ${\alpha}$ of $Y=\hat X\setminus X$ by boundedly saturated sets.
\endproclaim
We denote by $CPI'$ the condition $CPI$ without an assumption of metrizability
of $\hat X$. 
\proclaim{Theorem 9.6}
$CPII' \Leftrightarrow CPI'\Leftrightarrow CPI\Leftrightarrow HR\Rightarrow FW\Leftarrow
CPII$. 
\endproclaim
\demo{Proof}
By Proposition 9.1 the property of a compactification to be Higson dominated
is equivalent to the small action at infinity condition. Hence 
$HR\Leftrightarrow CPI'$. The exact sequence of pair 
$(\hat X,\hat X\setminus X)$ implies that $HR\Rightarrow FW$.
Clearly, $CPI\Rightarrow CPI'$. Show that
$CPI\Leftarrow CPI'$. Let $\hat X$ be a rationally acyclic compactification
of $X$. We can present $\hat X$ as the limit space of an inverse $\sigma$-continuous
system of metrizable compacta $\{\hat X_{\alpha}\mid \alpha\in A\}$. 
Like in the proof of
$2)\Rightarrow 3)$ of Lemma 9.3 we may assume that each $\hat X_{\alpha}$ is
an equivariant compactification of $X$. Consider the direct system
$\{H^i(\hat X_{\alpha};\Q)\mid \alpha\in A\}$. Let $M_{H^i}\subset A$ be the set defined in
the proof of the Dual Spectral Theorem. We recall that by the definition for any 
$\alpha\in M_{H^i}$, the kernel of $p_{\alpha}^*:H^i(\hat X_{\alpha};
\Q)\to H^i(\hat X;\Q)$ is trivial. Since $H^i(\hat X;\Q)=0$, it follows that
every $\hat X_{\alpha}$ is rationally acyclic.

Next we show that $CPII\Rightarrow FW$. Let $cX$ be 
$\min\{\hat X,\bar X=\nu X\cup X\}$ and let $g:\hat X\to cX$ be the natural
projection. First we show that every boundedly saturated open set 
$U\subset Y$ is of the form $g^{-1}(V)$. This implies that $V$ is open.
Show that
 $g(U)=g(Y)\setminus g(Y\setminus U)$. 
 Let $x\in g(U)$ and $x=g(z)$, $z\in U$.
Let $F$ be a closed set in $\hat X$ such that $F\cap(Y\setminus U)=\emptyset$ and
$z\in F$. Let $H$ be a closed set in $\hat X$ such that $H\cap Y=Y\setminus U$.
Note that if $g(F)\cap g(Y\setminus U)=\emptyset$, then 
$x\in g(Y)\setminus g(Y\setminus U)$. Assume the contrary
$g(F)\cap g(Y\setminus U)\ne\emptyset$. Then $g(F)\cap g(H)\ne\emptyset$
and hence  $(g(F)\cap X,g(H)\cap X)=g(F\cap X, H\cap X)\in\Cal P(cX)$.
Since $F\cap H=\emptyset$, we have $(F\cap X,H\cap X)\not\in\Cal P(\hat X)$.
Proposition 9.5 implies that $(F\cap X,H\cap X)\in\sP(\bar X)$. Then
there exist a constant $c$ and sequences $\{x_n\}\subset F\cap X$ and
$\{y_n\}\subset H\cap X$, tending to infinity, and such that $dist(x_n,y_n)<c$.
Since $U$ is boundedly saturated, we have 
$\overline{N_{2c}(F\cap X)}\cap Y\subset U$.
Hence $\overline{N_{2c}(F\cap X)}\cap (Y\setminus U)=\emptyset$. Hence
$(N_{2c}(F\cap X)\cap H)\setminus B_R(x_0)=\emptyset$ for large enough $R$.
This contradicts to the fact that $\{y_n\}$ tends to infinity.

Therefore for any covering $\alpha$ of $Y$ by boundedly saturated subsets
the image $g(\alpha)$ is an open cover of $cX\setminus X$ with the same nerve.
Consider the following diagram:
$$
\CD
H^{lf}_{k+1}(X;\Q) @= H^{lf}_{k+1}(X;\Q)\\
@V=VV  @V{\partial}VV\\
H_k(Y;\Q) @. H_k(cX\setminus X;\Q)\\
@V=VV  @VVV\\
H_k(\lim_{\leftarrow}N(\alpha);\Q) @>=>> H_k(\lim_{\leftarrow}N(g(\alpha));\Q)\\
\endCD
$$
The homomorphism $\partial$ is an equivariant split injection as a left
divisor of an equivariant isomorphism. For cohomology it is an
equivariant split surjection and hence $FW$ holds.

The condition $CPII'$ implies that the corresponding compactification is
Higson dominated.\qed
\enddemo

Thus the condition $FW$ is the weakest among above. The condition $FW$ is a
further coarsening of the following version of the Coarse Baum-Connes conjecture:
{\it the rational Roe index map $K^{lf}_*(X)\otimes\Q\to K_*(C(X))\otimes\Q$ is
an equivariant split monomorphism} which also implies the Novikov conjecture.

\head \S10 Open problems \endhead

{\bf Dimension and Higson corona.}
The connection between asymptotic dimensions of metric space and
 dimensions of its (nonmetrizable) Higson corona is not fully investigated
 yet. The problems in this paragraph are of particular interest for
 metric spaces of bounded geometry.
\proclaim{Problem 1}
Is it always true that $as\dim X=\dim\nu X$?
\endproclaim
In view of Theorem 7.2 and [D-K-U] this problem can be reformulated as
follows: {\it Does there exists a metric space $X$ of infinite asymptotic
dimension with $\dim \nu X<\infty$ ?} 
\proclaim{Problem 2}
Does there exists a metric space $X$ of an infinite asymptotic dimension
with $as\dim_*X<\infty$ ?
\endproclaim
A positive answer to Problem 2 gives a negative answer to Problem 1.
A relevant question is whether $as\dim_*X$ and $\dim\nu X$ always agree.
\proclaim{Problem 3}
Does the inequality $\dim_G\nu X\le as\dim_GX$ hold for all metric spaces $X$
and all abelian groups $G$ ?
\endproclaim
\proclaim{Problem 4}
Is it true that $\dim_{\Z}\nu\Gamma<\infty$ for all geometrically finite
groups $\Gamma$ ?
\endproclaim
An affirmative answer to Problem 3 would imply an affirmative answer to
Problem 4.

Using analogy one can define an asymptotic inductive dimension
$asInd$.
\proclaim{Problem 5}
What is the relation between $asInd X$ and other dimensions:
$as\dim X$, $as\dim_*X$ and $\dim\nu X$ ?
\endproclaim
\proclaim{Problem 6}
Let $X$ be a metric space of bounded geometry with slow dimension growth
and with $asdim X=\infty$.
Does it follow that $\dim\nu X=\infty$ ?
\endproclaim
\proclaim{Problem 7}
Can a metric space $X$ with $as\dim X=n$ be coarsely uniformly embedded
in a simply connected $(2n+1)$-dimensional non-positively curved manifold ?
\endproclaim
The answer is 'yes' for $2n+2$-dimensional manifolds.

\

{\bf Large scale Alexandroff Problem.}
\proclaim{Problem 8}
Does the equality $as\dim_{\Z}\Gamma=as\dim\Gamma$ holds for
geometrically finite groups $\Gamma$ ?
\endproclaim
The spectrum $\Bbb S=\{\Omega^{\infty}\Sigma^{\infty}S^n\}$ defines stable cohomotopy. Then
the stable cohomotopy dimension $\dim_{\Bbb S}X$ can be defined in terms of extensions of maps to
$\Omega^{\infty}\Sigma^{\infty}S^n$. Namely, $\dim_{\Bbb S}X\le n$ if and only if
every continuous map $f:A\to\Omega^{\infty}\Sigma^{\infty}S^n$  defined on a closed subset 
$A\subset X$ can be extended over  all $X$. It is easy to check that for every finite
$m$ extensions to $\Omega^m\Sigma^mS^n$ classify the covering dimension $\dim$.
Thus, $\dim X$ and $\dim_{\Bbb S}X$ are very close. One can define a macroscopic
version of the stable cohomological dimension as an asymptotic generalized cohomological
dimension.
\proclaim{Problem 9}
Does the asymptotic stable cohomotopy dimension $as\dim_{\Bbb S}\Gamma$
coincide with $as\dim\Gamma$ for  geometrically finite groups?
\endproclaim
Positive answers to Problems 8,9 imply the Novikov Conjecture.
\proclaim{Problem 10}
Does the equality $\dim_{\Z}X=\dim X$ hold for compact $H$-spaces ?
\endproclaim
\proclaim{Problem 11}
Does the equality $\dim_{\Z}G=\dim G$ hold for topological (noncompact)
groups $G$ ?
\endproclaim
\proclaim{Problem 12}
What it would be a coarse analogue of the Edwards resolution theorem?
\endproclaim

{\bf Absolute neighborhood extensors}.

\proclaim{Problem 13}
Is the space of probability measures $P(X)$ an absolute extensor in the
asymptotic topology?
\endproclaim
It is not difficult to show that $P(X)$ is AE in the class of finite dimensional
spaces or in the class of spaces with the slow dimension growth.
\proclaim{Problem 14}
Does the Homotopy Extension Theorem hold in the asymptotic category in full
generality?
\endproclaim
An affirmative answer to Problem 13 implies an affirmative answer to
Problem 14.
\proclaim{Problem 15}
Prove a macroscopic analog of the West theorem  stating that
every $ANE$ is homotopy eqivalent to a polyhedron.
\endproclaim

\

\

\

\

{\bf Fragments of the Micro-Macro topology dictionary.}

\ 

\ \ \ \ \ \  MICRO\ \ \ \ \ \ \ \ \ \ \ \ \ \ \ \ \ \ \ \ \ \ \ \ \ \ \ \ \ \ 
\ \ \ \ \ \ \ \ \ \ \ \ \ \ \ \ MACRO

1. Compactum\ \ \ \ \ \ \ \ \ \ \ \ \ \ \ \ \ \ \ \ \ \ \ 
\ \ \ \ \ \ \ \ \ \ \ \ \ Proper metric space of

i.e. compact metrizable space\ \ \ \ \ \ \ \ \ \ \ \ \ \ \ \ \ bounded geometry

\

2. Alexandroff-\v Cech approximation\ \ \ \ \ \ \ \ \ Anti-\v Cech approximation

by polyhedra\ \ \ \ \ \ \ \ \ \ \ \ \ \ \ \ \ \ \ \ \ \ \ \ \ \ \ \ 
\ \ \ \ \ \ \ \ \ \ \ by polyhedra

\

3. Lebesgue dimension $\dim$ \ \ \ \ \ \ \ \ \ \ \ \ \ \ \ \ \ \ 
\ Gromov dimensiona $as\dim$

\

4. Alexandroff characterization\ \ \ \ \ \ \ \ \ \ \ \ \ \ \ 
$\dim^c$=covering dimension of the

of $\dim$ by maps to $S^n$\ \ \ \ \ \ \ \ \ \ \ \ \ \ \ \ \ \ \ \ \ \ \
\ \ \ \ Higson corona

\

5. Local contractibility\ \ \ \ \ \ \ \ \ \ \ \ \ \ \ \ \ \ \ \ \ \ \ \ \ \
Uniform contractibility

\

6. Local $n$-connectivity\ \ \ \ \ \ \ \ \ \ \ \ \ \ \ \ \ \ \ \ \ \ \ \ \ 
Uniform $n$-connectivity

\

7. Neighborhood\ \ \ \ \ \ \ \ \ \ \ \ \ \ \ \ \ \ \ \ \ \ \ \ \ \ \ \ \ \ \ \ \ \
Asymptotic neighborhood

\

8. ANE\ \ \ \ \ \ \ \ \ \ \ \ \ \ \ \ \ \ \ \ \ \ \ \ \ \ \ \ \ \ \ \ \ \ 
\ \ \ \ \ \ \ \ \ \ \ \ ANE

\

9. Polyhedron\ \ \ \ \ \ \ \ \ \ \ \ \ \ \ \ \ \ \ \ \ \ \
\ \ \ \ \ \ \ \ \ \ \ \ \ \ Asymptotic polyhedron 

\ \ \ \ \ \ \ \ \ \ \ \ \ \ \ \ \ \ \ \ \ \ \ \ \ \ \ \ \ \ \ \ \ \ \ \ \
 \ \ \ \ \ \ \ \ \ \ \ \ \ \ \ \ \ \ \ \ \ \ \ \ \ \ \ \ or 
 
 \ \ \ \ \ \ \ \ \ \ \ \ \ \ \ \ \ \ \ \ \ \ \ \ \ \ \ \ \ \ \ \ \ \ \ \ \ 
 \ \ \ \ \ \ \ \ \ \ \ \ \ \ \ \ \ \ 
 open cone over a polyhedron
 
\

10. One point space\ \ \ \ \ \ \ \ \ \ \ \ \ \ \ \ \ \ \ \ \ \ \ \ \ \ \ \ \
\ 
$\R_+$ (as well as $\N$)

\

11. Unit interval $[0,1]$\ \ \ \ \ \ \ \ \ \ \ \ \ \ \ \ \ \ \ \ \ \ \ \ \ \ \ 
$\R^2_+$

\

12. $n$-Sphere $S^n$\ \ \ \ \ \ \ \ \ \ \ \ \ \ \ \ \ \ \ \ \ \ \ \ \ \ \ \ \ \ \ \
\ \ $\R^{n+1}$

\

13.\ \ \ ?\ \ \ \ \ \ \ \ \ \ \ \ \ \ \ \ \ \ \ \ \ \ \ \ \ \ \ \ \ \ \ \ 
\ \ \ \ \ \ \ \ \ \ \ \ \ \ \ Hyperbolic space $\H^n$

\

14. Homotopy\ \ \ \ \ \ \ \ \ \ \ \ \ \ \ \ \ \ \ \ \ \ \ \ \ \ \ \ \ \ \ \ 
\ \ \ \ \  Homotopy in $\sA$

\

15. Cohomology\ \ \ \ \ \ \ \ \ \ \ \ \ \ \ \ \ \ \ \ \ \ \ \ \ \ \ \ \ \ \ \ \ 
\ anti-\v Cech cohomology

\ \ \ \ \ \ \ \ \ \ \ \ \ \ \ \ \ \ \ \ \ \ \ \ \ \ \ \ \ \ \ \ \ \ \ \ \ \ 
\ \ \ \ \ \ \ \ \ \ \ \ \ \ \ \ \ \ Roe's cohomology

\

16. \v Cech homology\ \ \ \ \ \ \ \ \ \ \ \ \ \ \ \ \ \ \ \ \ \ \
\ \ \ \ \ \ \ Coarse homology

\ \ \ Steenrod homology

\

17. Cohomological dimension\ \ \ \ \ \ \ \ \ \ \ \ \ \ \ \ \ \  
Asymptotic cohomological dimension

\

18. Fundamental group\ \ \ \ \ \ \ \ \ \ \ \ \ \ \ \ \ \ \ \ \ \ \
Asymptotic fundamental 

\ \ \ \ \ \ \ \ \ \ \ \ \ \ \ \ \ \ \ \ \ \ \ \ \ \ \ \ \ \ \ \ \ \ \ \ \ 
\ \ \ \ \ \ \ \ \ \ \ \ \ \ \ \ \ group can be defined by using

\ \ \ \ \ \ \ \ \ \ \ \ \ \ \ \ \ \ \ \ \ \ \ \ \ \ \ \ \ \ \ \ \ \ \ \ \ 
\ \ \ \ \ \ \ \ \ \ \ \ \ \ \ \ \ $\R^2_+$ instead of $[0,1]$

\

19. Manifold\ \ \ \ \ \ \ \ \ \ \ \ \ \ \ \ \ \ \ \ \ \ \ \ \ \ \ \ \ \ \
\ \ \ \ \ \ Open contractible manifold?

\

20. \ \ \ \ ?\ \ \ \ \ \ \ \ \ \ \ \ \ \ \ \ \ \ \ \ \ \ \
\ \ \ \ \ \ \ \ \ \ \ \ \ \ \ \ \ \ \ \ \  Discrete group

\enddocument